\documentclass[MSNbibl,number,citesort,seceqn,dvips]{arxbj}
\usepackage{upgreek}
\usepackage{mathrsfs}

\aid{0}
\volume{19}
\issue{3}
\pubyear{2013}
\firstpage{954}
\lastpage{981}
\doi{10.3150/12-BEJ440} 

\makeatletter

\newremark{example}{Example}
\newremark{remark}{Remark}[section]
\newtheorem{lemma}{Lemma}[section]
\newtheorem{theorem}{Theorem}[section]
\newtheorem{corollary}{Corollary}[section]
\newremark{details}{Details}

\def\R{\mathbb{R}}
\def\x{\times}
\def\vt{\vartheta}
\def\ve{\varepsilon}
\def\el{\mathscr{R}_{n}}
\def\hel{\hat{\mathscr{R}}_{n}}
\def\1{\mathbf{1}}
\def\hz{\hat\zeta_n}

\newcommand{\trace}{\operatorname{trace}}
\newcommand{\sign}{\operatorname{sign}}
\newcommand{\Ex}{\operatorname{Ex}}

\def\Z{{\mathcal Z}}
\def\G{\mathbb G}
\def\F{\mathbb F}
\def
\prob{\mathscr{P}_n}
\newcommand{\implies}{\Longrightarrow}
\def\ldot{\dot\ell}
\def\sdot{\dot s}
\def\H{\mathbb H}

\makeatother

\begin{document}
\begin{frontmatter}

\title{Empirical likelihood approach to goodness of fit testing}
\runtitle{Goodness of Fit}

\begin{aug}
\author[1]{\fnms{Hanxiang} \snm{Peng}\thanksref{1}\ead[label=e1]{hpeng@math.iupui.edu}} \and
\author[2]{\fnms{Anton} \snm{Schick}\corref{}\thanksref{2}\ead[label=e2]{anton@math.binghamton.edu}}
\runauthor{H. Peng and A. Schick} 
\address[1]{Department of Mathematical Sciences, Indiana University
Purdue University at Indianapolis,
Indianapolis, IN 46202, USA. \printead{e1}}
\address[2]{Department of Mathematical
Sciences, Binghamton University, Binghamton, NY 13902, USA.\\ \printead{e2}}
\end{aug}

\received{\smonth{10} \syear{2010}}
\revised{\smonth{3} \syear{2012}}

%
\begin{abstract}
Motivated by applications to goodness of fit testing,
the empirical likelihood approach is generalized
to allow for the number of constraints
to grow with the sample size and for the constraints
to use estimated criteria functions.
The latter is needed to deal with nuisance parameters.
The proposed empirical likelihood based goodness of fit tests are
asymptotically distribution free. For univariate observations,
tests for a specified distribution, for a distribution of parametric form,
and for a symmetric distribution are presented.
For bivariate observations, tests for independence are developed.
\end{abstract}

%

\end{frontmatter}

\section{Introduction}

The empirical likelihood approach was introduced by Owen~\cite{O88,O90}
to construct
confidence intervals in a nonparametric setting, see also Owen~\cite{O01}.
As a likelihood approach possessing nonparametric properties, it does
not require us
to specify a distribution for the data and often yields more efficient
estimates of the parameters.
It allows data to decide the shape of confidence regions and is
Bartlett correctable (DiCiccio, Hall and Romano~\cite{DH91}).
The approach has been developed to various situations,
for example, to generalized linear models (Kolaczyk~\cite{K94}), local
linear smoother (Chen and Qin~\cite{CQ00}),
partially linear models (Shi and Lau~\cite{SL00}, Wang and Jing~\cite{WJ03}),
parametric and semiparametric models in multiresponse regression
(Chen and Van Keilegom~\cite{CK09}), linear regression with censored data
(Zhou and Li~\cite{ZL08}),
and plug-in estimates of nuisance parameters in estimating equations in
the context of
survival analysis (Qin and Jing~\cite{QJ01a}, Wang and Jing~\cite{WJ01}, Li and
Wang~\cite{LC01}).
Algorithms, calibration and higher-order precision of the approach can
be found
in Hall and La Scala~\cite{HS90}, Emerson and Owen~\cite{EO09} and Liu and Chen
\cite{LC10} among others.
It is especially convenient to incorporate side information expressed
through equality constraints.
Qin and Lawless~\cite{QL94} linked empirical likelihood with finitely many
estimating equations.
These estimating equations serve as finitely many equality constraints.

In semiparametric settings, information on the model can often be expressed
by means of \emph{infinitely} many constraints which may also depend on
parameters of the model. In goodness of fit testing, the null
hypothesis can
typically be expressed by infinitely many such constraints.
This is the case
when testing for a fixed distribution (see Example~\ref{ex.1} below),
when testing for a given parametric model (Example~\ref{ex.model}),
when testing for symmetry about a fixed point (Example~\ref{ex.symm}),
and when testing for independence (Example~\ref{ex.indep}).
Modeling conditional expectations
can also be done by means of infinitely many constraints.
This has applications to heteroscedastic regression models (Section
\ref{sim}) and to
conditional moment restriction models treated by Tripathi and Kitamura~\cite{TK03}
using a smoothed empirical likelihood approach.

Recently Hjort, McKeague and Van Keilegom~\cite{HMV09} extended the scope of
the empirical method. In particular,
they developed a general theory for constraints with nuisance
parameters and
considered the case with infinitely many constraints.
Their results for infinitely many constraints, however, do not allow
for nuisance parameters.
In this paper we will fill this gap and in the process improve on their results.
Let us now discuss some of our results in the following special case.

Let $Z_1,\dots,Z_n$ be independent copies of a random vector $Z$ with
distribution $Q$. Let $u_1,u_2,\dots$ be orthonormal elements of
\[
L_{2,0}(Q)= \biggl\{ u\in L_2(Q)\dvt \int u \,
\mathrm{d}Q=0 \biggr\}.
\]
Then the random variables $u_1(Z),u_2(Z),\dots$ have mean zero,
variance one and
are uncorrelated.
Now consider the empirical likelihood based on the first $m$ of these functions,
\[
\el= \sup \Biggl\{\prod_{j=1}^n n
\pi_j\dvt\pi\in\prob, \sum_{j=1}^n
\pi_ju_k(Z_j)=0, k=1,\dots,m \Biggr\},
\]
where $\prob= \{ \pi= (\pi_1,\dots,\pi_n)^{\top}\in[0,1]^n\dvt  \pi_1+\cdots+\pi_n=1\}$
denotes the closed probability simplex in dimension $n$.
For fixed $m$, it follows from Owen's work that $-2\log\el$ has asymptotically
a chi-square distribution with $m$ degrees of freedom. In other words,
%
\begin{equation}
\label{wilk} P \bigl( -2\log\el> \chi^2_{1-\alpha}(m) \bigr) \to
\alpha,\qquad 0 < \alpha<1,
\end{equation}
where $\chi^2_{\beta}(m)$ denotes the $\beta$-quantile of the chi-square
distribution with $m$ degrees of freedom.
Hjort \textit{et al.}~\cite{HMV09} have shown that (\ref{wilk}) holds under some
additional assumptions
even if $m$ tends to infinity with $n$ by proving the
asymptotic normality result
%
\begin{equation}
\label{hmv} 
(-2 \log\el- m)/{\sqrt{2m}}\quad \implies \quad N(0,1).
\end{equation}
This result requires higher moment assumptions on the functions
$u_1,u_2,\dots$ and restrictions on the rate at which $m$
can tend to infinity. For example, if the functions $u_1,u_2,\dots$
are uniformly bounded, then the rate $m^3=\mathrm{o}(n)$ suffices for (\ref{hmv}).
They also state in their Theorem 4.1, that if $\sup_{k} \int|u_k|^q\, \mathrm{d}Q$
is finite for some $q>2$, then $m^{3+6/(q-2)}=\mathrm{o}(n)$ suffices for (\ref{hmv}).
A gap in their argument was fixed by Peng and Schick~\cite{PS12}.
We shall show that larger $m$ are allowed in some cases.
In particular, for $q=4$, it suffices that $m^4=\mathrm{o}(n)$ holds
(instead of their $m^6=\mathrm{o}(n)$) and if $q=3$, then $m_n^6=\mathrm{o}(n)$ is enough
(instead of their $m_n^9=\mathrm{o}(n)$), see our Theorems~\ref{thm.3b}
and~\ref{thm.3d} below.

Our rate $m^4=\mathrm{o}(n)$ for $q=4$ matches the rate given in Theorem 2
of Chen, Peng and Qin~\cite{CPQ09}. These authors obtain asymptotic normality
for $m$ larger than in Hjort \textit{et al.}~\cite{HMV09} by imposing additional
structural assumptions.
These assumptions, however, are typically not met in the applications
we have in mind.

One of the key points in our proof is a simple condition for
the convex hull of some vectors $x_1,\dots,x_n$ to have the origin
as an interior point. Our condition is that the smallest eigenvalue of
$\sum_{i=1}^n x_ix_i^{\top}$ exceeds $5|\sum_{i=1}^n x_i|\max_{1\leq j \leq n}|x_j|$.
Here, $|x|$ denotes the euclidean norm of a vector $x$.
This sufficient condition ties in nicely with the other requirements
used to establish the asymptotic behavior of the empirical likelihood
and is typically implied by these. For example,
conditions (A1)--(A3) in Theorem 2.1 of Hjort \textit{et al.}~\cite{HMV09} already imply
their (A0).
Thus, the conclusion of their theorem is valid under (A1)--(A3) only,
see our Theorem~\ref{thm.1}.

Let us now look at the case when the functions $u_1,u_2,\dots$ are unknown.
Then we can work with the empirical likelihood
\[
\hel= \sup \Biggl\{\prod_{j=1}^n n
\pi_j\dvt\pi\in\prob, \sum_{j=1}^n
\pi_j\hat u_k(Z_j)=0, k=1,\dots,m \Biggr
\},
\]
where $\hat u_k$ is an estimator of $u_k$ such that
%
\begin{equation}
\label{p0} \sum_{k=1}^m \frac{1}{n}
\sum_{j=1}^n \bigl(\hat
u_k(Z_j)-u_k(Z_j)
\bigr)^2 = \mathrm{o}_p \bigl(m^{-1} \bigr).
\end{equation}
Now, we have the conclusion
$(-2\log\hel- m)/\sqrt{2m} \implies N(0,1)$
under the condition
%
\begin{equation}
\label{p1} \sum_{k=1}^m \Biggl(
n^{-1/2}\sum_{j=1}^n \bigl(\hat
u_k(Z_j)-u_k(Z_j) \bigr)
\Biggr)^2 = \mathrm{o}_p(1)
\end{equation}
and mild additional conditions such as
\begin{enumerate}[(ii)]
\item[(i)] $|\hat u_k|+|u_k| \le B$ for some constant $B$ and all $k$ and $m^3=\mathrm{o}(n)$,
or
\item[(ii)] $\sum_{k=1}^m \int u_k^4\, \mathrm{d}Q = \mathrm{O}(m^2)$ and $m_n^4=\mathrm{o}(n)$.
\end{enumerate}
Our results, however, go beyond this simple result.
If (\ref{p1}) is replaced by
%
\begin{equation}
\label{p2} \sum_{k=1}^m \Biggl(
n^{-1/2}\sum_{j=1}^n \bigl(\hat
u_k(Z_j)-u_k(Z_j)+ E
\bigl[u_k(Z)\psi^{\top
}(Z) \bigr]\psi(Z_j)
\bigr) \Biggr)^2 = \mathrm{o}_p(1)
\end{equation}
with $\psi$ a measurable function into $\R^q$ which is standardized
under $Q$ in the sense that $E[\psi(Z)]=0$ and $E[\psi(Z)\psi^{\top
}(Z)]=I_q$,
the $q\times q$ identity matrix,
then the conclusion
$(-2\log\hel- (m-q))/\sqrt{2(m-q)} \implies N(0,1)$ holds
under (i) or (ii).

Our paper is organized as follows.
In Section~\ref{ex}, we give four examples that motivate our research.
The emphasis in these examples is on goodness of fit testing.
The proposed empirical likelihood based goodness of fit tests
are asymptotically distribution free.
For univariate observations, tests for a specified distribution,
for a distribution of parametric form,
and for a symmetric distribution are presented.
For bivariate observations, tests for independence are discussed.
Another example is given in Section~\ref{sim} with a small simulation study.
This example considers tests for the regression parameters
in simple linear heteroscedastic regression.
The simulations compare our new procedure based on infinitely many
constraints with the classical empirical likelihood procedure
and illustrate improvements by the new procedures.
In Section~\ref{not}, we introduce notation and recall some results on
the spectral norm of matrices.
In Section~\ref{max}, we derive a lemma that extracts the essence
from the proofs of Owen (\cite{O01}, Chapter 11)
and also obtains the
aforementioned sufficient condition for a convex hull of vectors
to contain the origin as interior point.
The results are derived for non-stochastic vectors
and formulated as inequalities.
The inequalities are used in Section~\ref{app} to obtain the behavior
of the
empirical likelihood with random vectors whose dimension may increase.
The results are formulated abstractly and do not require independence.
In Section~\ref{main}, we specialize our results to the case of independent
observations with infinitely many constraints, both known and unknown.
We also briefly discuss the behavior under contiguous alternatives.
The details for our examples are given in Section~\ref{det}.

\section{Motivating examples} \label{ex}

In this section, we give examples that motivated the research in this paper.

\begin{example}[(Testing for a fixed distribution)]\label{ex.1}
Let $X_1,\dots,X_n$ be independent copies of a random variable $X$.
Suppose we want to test whether their common distribution function $F$
equals a known \emph{continuous} distribution function $F_0$.
Under the null hypothesis, we have $E[h(X)]=0$ for every $h\in L_{2,0}(F_0)$,
and $F_0(X)$ has a uniform distribution.
An orthonormal basis of $L_{2,0}(F_0)$ is thus given by
$v_1\circ F_0,v_2\circ F_0,\dots$ for any orthonormal basis
$v_1,v_2,\dots$ of $L_{2,0}(U)$,
where $U$ is the uniform distribution on $(0,1)$.
We shall work with the trigonometric basis $\phi_1,\phi_2,\dots$
defined by
%
\begin{equation}
\label{trig-basis} \phi_k(x)= \sqrt{2} \cos(k\uppi x),\qquad x
\in[0,1], k=1,2,\dots,
\end{equation}
as these basis functions are uniformly bounded by $\sqrt{2}$.
As test statistic, we take
\[
\el(F_0)= \sup \Biggl\{\prod
_{j=1}^n n\pi_j\dvt\pi\in\prob, \sum
_{j=1}^n \pi_j
\phi_k \bigl(F_0(X_j) \bigr)=0, k=1,\dots,m
\Biggr\}
\]
which uses the first $m$ of the trigonometric functions.
Under the null hypothesis, we have
$
P(-2\log\el(F_0) > \chi^2_{1-\alpha}(m)) \to\alpha
$
for every $0< \alpha<1$ as both $m$ and $n$ tend to infinity
and $m^3/n$ tends to zero.
Thus, the test $\1[-2\log\el(F_0) > \chi^2_{1-\alpha}(m)]$
has asymptotic size~$\alpha$.
Here, we are still in the framework of Hjort \textit{et al.}~\cite{HMV09} with
infinitely many \emph{known} constraints.
\end{example}

\begin{example}[(Testing for a parametric model)]\label{ex.model}
Let $X_1,\dots,X_n$ be again independent and identically distributed
random variables.
But now suppose we want to test whether their common distribution
function $F$ belongs to a
model $\mathscr{F}=\{ F_{\vt}\dvt  \vt\in\Theta\}$ indexed by an open
subset $\Theta$ of $\R^q$.
Suppose that the distribution functions $F_{\vt}$ have densities
$f_{\vt}$ such that
the map $\vt\mapsto s_{\vt}=\sqrt{f_{\vt}}$ is continuously
differentiable in $L_2$
with derivative $\vt\mapsto\sdot_{\vt}$ and the matrix
$J(\vt)= 4 \int\sdot_{\vt}(x)\sdot_{\vt}(x)^{\top} \, \mathrm{d}x$
is invertible for each $\vt\in\Theta$.
In this case we set $\ldot_{\vt}=2\sdot_{\vt}/s_{\vt}$.
Let now $\hat\theta$ be an estimator of the parameter in the model.
We require it to satisfy the stochastic expansion
%
\begin{equation}
\label{para-eff} \hat\theta= \theta+ \frac{1}{n}\sum
_{j=1}^nJ(\theta)^{-1}
\ldot_{\theta}(X_j) + \mathrm{o}_{P_{\theta}}
\bigl(n^{-1/2} \bigr)
\end{equation}
for each $\theta\in\Theta$, where $P_{\theta}$ is the measure for
which $F=F_{\theta}$.
Such estimators are efficient in the parametric model.
Candidates are maximum likelihood estimators.
As test statistic we take $\el(F_{\hat\theta})$,
the test statistic from the previous example with $F_0$ replaced by
$F_{\hat\theta}$.
Here, we are no longer in the framework of Hjort \textit{et al.}~\cite{HMV09} as we now have
infinitely many
\emph{unknown} constraints.
We shall show that under the null hypothesis
$
P(-2\log\el(F_{\hat\theta}) > \chi^2_{1-\alpha}(m-q)) \to\alpha
$
for every $0 < \alpha< 1$
as both $m$ and $n$ tend to infinity and $\log n m^3/n$ tends to zero.
In view of this result, the test
$\1[-2\log\el(F_{\hat\theta}) > \chi^2_{1-\alpha}(m-q)]$
has asymptotic size $\alpha$.
It is crucial for our result that we have chosen an estimator $\hat
\theta$
satisfying~(\ref{para-eff}).
\end{example}

\begin{example}[(Testing for symmetry)]\label{ex.symm}
Let $X_1,\dots,X_n$ be independent copies of a random variable $X$
with a continuous distribution function $F$. We want to test whether $F$
is symmetric about zero in the sense that $F(t)=1-F(-t)$ for all real $t$.
Under the null hypothesis of symmetry, the random variables $\sign(X)$ and
$|X|$ are independent, and $\sign(X)$ takes values $-1$ and $1$
with probability one half.
This is equivalent to $E[\sign(X) v(|X|)]=0$ for every
$v\in L_2(H)$, where $H$ is the distribution function of $|X|$.
Since $H$ is continuous, an orthonormal system of $L_2(H)$ is given by
$\phi_0\circ H, \phi_1\circ H, \dots$ where $\phi_0=1$ and
$\phi_1,\phi_2,\dots$ are given in (\ref{trig-basis}).
This suggests the test statistic
\[
\begin{aligned} \el= \sup \Biggl\{\prod
_{j=1}^n n\pi_j\dvt\pi\in\prob, \sum
_{j=1}^n \pi_j
\sign(X_j)\phi_k(R_j)=0, k=0,\dots,m \Biggr
\}, \end{aligned} %
\]
where $R_j=\H(|X_j|)$ and $\H$ is the empirical distribution function based
on $|X_1|,\dots,|X_n|$.
We shall show that under symmetry one has
$
P(-2 \log\el> \chi^2_{1-\alpha}(m+1)) \to\alpha
$
for every $0<\alpha<1$
as $m$ and $n$ tend to infinity and $m^3/n$ tends to zero.
From this, we derive that the test $\1[-2\log\el> \chi^2_{1-\alpha}(m+1)]$
has asymptotic size $\alpha$.
\end{example}

\begin{example}[(Testing for independence)]\label{ex.indep}
Let $(X_1,Y_1),\dots,(X_n,Y_n)$ be independent copies of a
bivariate random vector $(X,Y)$.
We assume that the marginal distribution functions
$F$ and $G$ are continuous.
We want to test whether $X$ and $Y$ are independent.
Independence is equivalent to $E[a(X)b(Y)]=0$
for all $a\in L_{2,0}(F)$ and $b\in L_{2,0}(G)$
and thus equivalent to
$E[\phi_k(F(X))\phi_l(G(Y))]=0$ for all positive integers $k$ and $l$.

\begin{longlist}[(b)]
\item[(a)] Assume first that $F$ and $G$ are known.
This is for example the case in an actuarial setting
where $X$ and $Y$ denote residual lifetimes
and their distribution functions are available from life tables.
Motivated by the above, we take as test statistics
\[
\el(F,G) = \sup \Biggl\{\prod_{j=1}^n n
\pi_j\dvt\pi\in\prob, \sum_{j=1}^n
\pi_j\phi_k \bigl(F(X_j) \bigr)
\phi_l \bigl(G(Y_j) \bigr)=0, k,l=1,\dots,r \Biggr\}.
\]
Under the null hypothesis, one has
$
P(-2 \log\el(F,G) > \chi^2_{1-\alpha}(r^2)) \to\alpha
$
for every $0<\alpha<1$ as $r$ and $n$ tend to infinity and $r^6/n$
tends to zero.
Here, we are in the framework of Hjort, McKeague and Van Keilegom~\cite{HMV09}.
The above shows that the test $\1[-2\log\el(F,G) > \chi^2_{1-\alpha}(r^2)]$
has asymptotic size $\alpha$.

\item[(b)] Now assume that $F$ and $G$ are unknown.
In this case, we replace both marginal distribution functions by their
empirical distribution functions. The resulting test statistic is
$\el(\F,\G)$, where $\F$ denotes the empirical distribution based
on $X_1,\dots,X_n$ and $\G$ the one based on $Y_1,\dots,Y_n$.
We shall show that under the null hypothesis
$
P(-2 \log\el(\F,\G) > \chi^2_{1-\alpha}(r^2)) \to\alpha
$
for every $0<\alpha<1$ as $r$ and $n$ tend to infinity and $r^6/n$
tends to zero.
Thus the test $\1[-2\log\el(\F,\G) > \chi^2_{1-\alpha}(r^2)]$
has asymptotic size $\alpha$.
\end{longlist}
\end{example}

\begin{remark}\label{rem.reg}
Suppose that $(X,Y)$ form a simple linear homoscedastic regression model,
$Y= \beta_1+\beta_2 X +\ve$, with $X$ and $\ve$ independent.
We can use the test statistic from case (b) to test the hypothesis
whether the slope parameter $\beta_2$ is zero.
Indeed, $\beta_2=0$ is equivalent to the independence of $X$ and $Y$.
\end{remark}

\begin{remark}\label{rem.alter}
The asymptotic distributions of the above tests under contiguous
alternatives are linked to non-central chi-square distributions;
see Remark~\ref{rem.alt} for details. As the non-centrality parameters
are bounded,
the local asymptotic power along such a contiguous alternative
coincides with the level. Our tests are asymptotically equivalent
to Neyman's smooth tests~\cite{N39} with increasing dimensions.
In view of the optimality results of Inglot and Ledwina~\cite{IL96},
for those tests under moderate deviations, we
expect similar results for our tests. Of course, this needs to be explored
more carefully.
\end{remark}

\section{Another example and simulations}\label{sim}

Let $(X_1,Y_1),\dots,(X_n,Y_n)$ be independent copies of $(X,Y)$, where
$Y=\beta_1+ \beta_2 X + \ve$, with $E[\ve|X]=0$, $\sigma^2(X)=E[\ve^2|X]$
bounded and bounded away from zero, and $E[\ve^4]<\infty$.
Assume that $X$ has a finite variance and
a continuous distribution function $G$.
We are interested in testing whether the regression parameter
$\beta=(\beta_1,\beta_2)^{\top}$ equals some specific value $\theta$.
We could proceed as in Owen~\cite{O91} and use the test
$
\delta_0 = \1[-2 \log{\mathscr{R}}_{n0}(\theta) > \chi^2_{1-\alpha}(2)]
$
based on the empirical likelihood
\[
\mathscr{R}_{n0}(\theta)= \sup \Biggl\{\prod
_{j=1}^n n\pi_j\dvt\pi\in \prob, \sum
_{j=1}^n \pi_j %
\pmatrix{ 1
\cr
X_j } %
(Y_j-
\theta_{1}-\theta_{2} X_j)=0 \Biggr\}.
\]
But this empirical likelihood does not use all the information of the model.
Here we have $E[a(X)\ve]=0$ for every $a\in L_2(G)$.
Since $G$ is continuous (but unknown), we work with the empirical likelihood
\[
\hat{\mathscr{R}}_{n1}(\theta) = \sup \Biggl\{\prod
_{j=1}^n n\pi_j\dvt \pi\in\prob, \sum
_{j=1}^n \pi_ju_r
\bigl(\G(X_j) \bigr) (Y_j-\theta_{1}-
\theta_{2} X_j)=0 \Biggr\},
\]
where $u_r=(1,\phi_1,\dots,\phi_r)^{\top}$ and
$\G$ is the empirical distribution function based on the
covariate observations $X_1,\dots,X_n$.
It follows from Corollary~\ref{cor.5} and Lemma~\ref{lm.tech} below that
$P(-2 \log\hat{\mathscr{R}}_{n1}(\theta)>\chi^2_{1-\alpha}(1+r))
\to\alpha$
if $r^4 =\mathrm{o}(n)$.
The resulting test is
$
\delta_1 = \1[-2 \log\hat{\mathscr{R}}_{n1}(\theta) >\chi^2_{1-\alpha}(r+1)].
$
Both tests have asymptotic size $\alpha$.

We performed a small simulation study to compare the procedures.
For our simulation, we chose $\alpha=0.05$ and $n=100$ and took
$\theta=(1,2)^{\top}$.
We modeled the error $\varepsilon$ as $s(X)\eta$,
with $s(X)= \min(\sqrt{1+X^2},100)$ and $\eta$ independent of $X$.
As distributions for $X$, we chose the exponential distribution
with mean 5 ($\Ex$(5)) and the $t$-distribution with three
degrees of freedom ($t$(3)), while for $\eta$
we chose the standard normal distribution ($N$(0,~1)) and
the double exponential distribution with location 0 and scale 0.5 ($L$(0, 0.5)).
%
\begin{table}
\tabcolsep=0pt
\caption{Simulated powers of the tests $\delta_0$ and
$\delta_1$}\label{table:power-i-a}
\begin{tabular*}{\textwidth}{@{\extracolsep{\fill}}lllllllllllll@{}}
\hline
& & &\multicolumn{5}{l}{$t(3)$} & \multicolumn{5}{l@{}}{$\Ex(5)$} \\ [-5pt]
& & &\multicolumn{5}{c}{\hrulefill} & \multicolumn{5}{c@{}}{\hrulefill} \\
&\multicolumn{1}{l}{$\beta_1$} &
\multicolumn{1}{l}{$\beta_2$} &
\multicolumn{1}{l}{0} &
\multicolumn{1}{l}{2} &
\multicolumn{1}{l}{3} &
\multicolumn{1}{l}{4} &
\multicolumn{1}{l}{5} &
\multicolumn{1}{l}{0} &
\multicolumn{1}{l}{2} &
\multicolumn{1}{l}{3} &
\multicolumn{1}{l}{4} &
\multicolumn{1}{l@{}}{5} \\
\hline
$N(0, 1)$& 0.6 & 2.3
& 0.71 
& 0.88 & 0.86 & 0.85 & 0.84
& 0.38 
& 0.37 & 0.39 & 0.40 & 0.41 \\
& 0.8 & 1.5
& 0.68 
& 0.82 & 0.84 & 0.83 & 0.83
& 0.95 
& 0.99 & 0.99 & 0.99 & 0.99 \\
& 1.0 & 2.0
& 0.13 
& 0.09 & 0.10 & 0.12 & 0.13
& 0.12 
& 0.07 & 0.09 & 0.12 & 0.14 \\
& 1.2 & 2.2
& 0.37 
& 0.42 & 0.43 & 0.43 & 0.44
& 0.51 
& 0.54 & 0.52 & 0.50 & 0.52\\
& 1.4 & 1.7
& 0.71 
& 0.88 & 0.87 & 0.86 & 0.86
& 0.37 
& 0.34 & 0.37 & 0.40 & 0.44 \\
[6pt]
$L(0, 0.5)$& 0.6 & 2.3
& 0.89 
& 0.98 & 0.99 & 0.98 & 0.98
& 0.61 
& 0.64 & 0.68 & 0.71 & 0.74\\
& 0.8 & 1.5
& 0.84 
& 0.96 & 0.98 & 0.98 & 0.98
& 0.93 
& 1.00 & 1.00 & 1.00 & 1.00 \\
& 1.0 & 2.0
& 0.14 
& 0.10 & 0.14 & 0.17 & 0.21
& 0.13 
& 0.10 & 0.11 & 0.14 & 0.17\\
& 1.2 & 2.2
& 0.57 
& 0.70 & 0.70 & 0.70 & 0.74
& 0.68 
& 0.84 & 0.84 & 0.82 & 0.83 \\
& 1.4 & 1.7
& 0.89 
& 0.99 & 0.99& 0.99 & 0.99
& 0.62 
& 0.67 & 0.72 & 0.73 & 0.76\\
\hline
\end{tabular*}
\end{table}

Table~\ref{table:power-i-a} reports simulated powers of the tests $\delta_0$ and $\delta_1$
(with several choices of $r$) and for some values of $\theta$.
The reported values are based on 1000 repetitions.
The column labeled 0 corresponds to Owen's test $\delta_0$,
while the columns labeled 2, 3, 4, 5 correspond to our tests $\delta_1$
with $r=2,3,4,5$, respectively.
Clearly our new test is more powerful than the traditional test.
The values in the rows corresponding to the parameter values
$(1.0, 2.0)$ are the observed significance levels of the nominal
significance level $0.05$. Our new test overall has closer observed
significance levels than the traditional one except for $r=5$.

\section{Notation} \label{not}
In this section, we introduce some of the notation we use throughout.
We write $|A|$ for the euclidean norm and $|A|_o$ for the operator
(or spectral) norm of a matrix $A$ which are defined by
\[
|A|^2 = \trace \bigl(A^{\top} A \bigr) = \sum
_{i,j} A_{ij}^2 \quad\mbox{and}\quad
|A|_o = \sup_{|u|=1} |Au| = \sup_{|u|=1}
\bigl(u^{\top} A^{\top}Au \bigr)^{1/2}.
\]
In other words, the squared euclidean norm $|A|^2$ equals the sum
of the eigenvalues of $A^{\top}A$,
while the squared operator norm $|A|_o^2$ equals
the largest eigenvalue of $A^{\top}A$.
Consequently, the inequality $|A|_o \le|A|$ holds. Thus, we have
\[
|Ax|\leq|A|_o |x| \leq|A||x|
\]
for compatible vectors $x$. We should also point out the identity
\[
|A|_o= \sup_{|u|=1}\sup_{|v|=1} u^{\top} A v.
\]
If $A$ is a nonnegative definite symmetric matrix, this simplifies to
\[
|A|_o= \sup_{|u|=1} u^{\top}Au.
\]
Using this and the Cauchy--Schwarz inequality, we obtain
\begin{eqnarray}
\label{onorm.csi} \biggl|\int fg^{\top}\, \mathrm{d}\mu \biggr|^2_o
&\leq& \biggl|\int ff^{\top}\, \mathrm{d}\mu \biggr|_o \biggl|\int
gg^{\top} \, \mathrm{d}\mu \biggr|_o,
\\
\label{onorm.cov} \biggl|\int ff^{\top}\, \mathrm{d}\mu
\biggr|_o &\leq& \int|f|^2 \,\mathrm{d}\mu,
\end{eqnarray}
whenever $\mu$ is a measure and $f$ and $g$ are measurable functions into
$\R^{s}$ and $\R^t$ such that $\int|f|^2 \,\mathrm{d}\mu$
and $\int|g|^2 \,\mathrm{d}\mu$ are finite.
As a special case, we derive the inequality
\[
|S_{x+y}-S_x|_o \leq|S_y|_o
+ 2 |S_x|_o^{1/2} |S_y|_o^{1/2}
\]
and therefore
%
\begin{equation}
\label{onorm.diff} 
|S_{x+y}-S_x|_o
\leq\frac{1}{n}\sum_{j=1}^n|y_i|^2
+ 2|S_x|_o^{1/2} \Biggl(\frac{1}{n}\sum
_{j=1}^n|y_i|^2
\Biggr)^{1/2} 
\end{equation}
with
\[
S_{x+y}=\frac{1}{n}\sum_{j=1}^n(x_j+y_j)
(x_j+y_j)^{\top},\qquad S_{x}=
\frac{1}{n}\sum_{j=1}^nx_jx_j^{\top},
\qquad S_y=\frac{1}{n}\sum_{j=1}^ny_jy_j^{\top}
\]
for vectors $x_1,y_1,\dots,x_n,y_n$ of the same dimension.

\section{A maximization problem} \label{max}
Let $x_1, \ldots, x_n$ be $m$-dimensional vectors. Set
$x_*=\max_{1\leq j \leq n} |x_j|$,
\[
\bar x=\frac{1}{n}\sum_{j=1}^nx_j,
\qquad 
S=\frac{1}{n}\sum_{j=1}^nx_j
x_j^\top,\qquad x^{(\nu)}=\sup_{|u|=1} \Biggl|
\frac{1}{n}\sum_{j=1}^n
\bigl(u^{\top
}x_j \bigr)^{\nu} \Biggr|,\qquad \nu=3,4,
\]
and let $\lambda$ and $\Lambda$ denote the smallest and largest eigenvalue
of the matrix $S$,
\[
\lambda= \inf_{|u|=1} u^{\top} S u\quad \mbox{and}\quad \Lambda=
\sup_{|u|=1} u^{\top} S u.
\]
%
Using Lagrange multipliers, Owen~\cite{O88,O01} obtained the identity
\[
\mathscr{R} = \sup \Biggl\{\prod_{j=1}^n
n \pi_j\dvt\pi\in\prob, \sum_{j=1}^n
\pi_jx_j=0 \Biggr\}= \prod
_{j=1}^n \frac{1}{1+\zeta^{\top}x_j}
\]
if there exists a $\zeta$ in $\R^m$ such that
$1+\zeta^{\top}x_j>0$, $j=1,\dots,n$,
and
%
\begin{equation}
\label{o2} \sum_{j=1}^n
\frac{x_j}{1+\zeta^{\top} x_j}=0.
\end{equation}
He also showed that such a vector $\zeta$ exists and is unique if
(i) the origin is an interior point of the convex hull of
$x_1,\dots,x_n$ and (ii) the matrix $S$ is invertible.
Let us now show that the inequality $\lambda> 5x_*|\bar x|$ implies
these two conditions. Indeed, the matrix $S$ is then positive definite
and hence invertible as its smallest eigenvalue $\lambda$ is positive.
To show (i), we will rely on the following lemma.

\begin{lemma}\label{lm.lbd}
A random variable $Y$ with $E[Y]=0$ and $P(|Y|\le K)=1$ for some
positive $K$ obeys the inequality
\[
P(Y>a) \ge\frac{E[Y^2]-2Ka}{2K^2}, \qquad 0\leq a <K.
\]
\end{lemma}

\begin{pf}
Fix $a$ in $[0,K)$. By the properties of $Y$,
we obtain
$
2K^2 P(Y>a)\geq2KE[Y\1[Y>a]] \geq2KE[Y\1[Y>0]]-2Ka
$
and $2KE[Y\1[Y>0]] = KE[|Y|] \ge E[Y^2]$.
\end{pf}

The origin is an interior point of the
convex hull of $x_1,\dots, x_n$ if for every unit vector
$u \in\R^m$ there is at least one $j\in\{1,\dots, n\}$
such that $u^\top x_j>0$.
This latter condition is equivalent to
\[
N= \inf_{|u|=1} \sum_{j=1}^n \1
\bigl[u^\top x_j >0 \bigr]\geq1.
\]
For a unit vector $u$, we have $-u^{\top} \bar x \leq|\bar x|$ and
thus
\[
\sum_{j=1}^n \1 \bigl[u^\top
x_j >0 \bigr] \geq \sum_{j=1}^n
\1 \bigl[u^\top(x_j-\bar x)> |\bar x| \bigr] = N(u).
\]
It follows from the triangle inequality that
$|x_j -\bar x| \leq|x_j|+|\bar x|\leq2x_*$ for $j=1,\dots,n$.
Note that $x_*$ is positive if $S$ is positive definite.
Thus, Lemma~\ref{lm.lbd} yields the lower bound
$N(u)/n \ge(\sigma^2(u)-4x_* |\bar x|)/(8x_*^2)$ with
%
\[
\sigma^2(u) =\frac{1}{n}\sum_{j=1}^n
\bigl(u^{\top}(x_j-\bar x) \bigr)^2 =
u^T Su - \bigl(u^{\top} \bar x \bigr)^2 \geq
\lambda- |\bar x|^2 \geq\lambda- x_*|\bar x|.
\]
Thus, we have $N\ge n(\lambda- 5|\bar x|x_*)/(8x_*^2)$.
This shows that the inequality $\lambda> 5|\bar x| x_*$
implies $N\geq1$ and hence the desired condition (i).

Assume now that the inequality $\lambda> 5 x_*|\bar x|$ holds.
We proceed as on page 220 of Owen~\cite{O01}. Let $u$ be a unit vector such that
$\zeta=|\zeta|u$. Then we have the identity
\[
0= \frac{1}{n}\sum_{j=1}^n
\frac{u^{\top}x_j (1+\zeta^{\top}x_j
-\zeta^{\top
}x_j)}{1+\zeta^{\top} x_j} = u^{\top}\bar x - |\zeta| \frac{1}{n}\sum
_{j=1}^n\frac{(u^{\top
}x_j)^2}{1+\zeta^{\top}x_j}
\]
and the inequality
\[
\lambda\leq u^{\top}S u = \frac{1}{n}\sum
_{j=1}^n \bigl(u^{\top}x_j
\bigr)^2 \leq\frac{1}{n}\sum_{j=1}^n
\frac{(u^{\top}x_j)^2 (1+|\zeta|
x_*)}{1+\zeta^{\top}x_j}.
\]
Consequently, we find
$
\lambda|\zeta| \leq(1+|\zeta| x_*) u^{\top}\bar x \leq(1+|\zeta|
x_*) |\bar x|
$
and obtain the bound
%
\begin{equation}
\label{o3} |\zeta|\leq\frac{|\bar x|}{\lambda-|\bar x| x_*}.
\end{equation}
From this, one immediately derives
\begin{eqnarray}
\label{o4} |\zeta|x_* &\leq&\frac{|\bar x|x_*}{\lambda-|\bar x| x_*} < \frac{1}{4},
\\
\label{o5} \max_{1\leq j \leq  n}\frac{1}{1+\zeta^{\top}x_j} &\leq&\frac{1}{1-|\zeta|x_*} <
\frac{4}{3},
\\
\label{o6} \frac{1}{n}\sum_{j=1}^n
\bigl(\zeta^{\top}x_j \bigr)^2 &=&
\zeta^{\top} S \zeta \le\Lambda|\zeta|^2 \leq
\frac{\Lambda|\bar x|^2}{(\lambda- |\bar x|x_*)^2}.
\end{eqnarray}
The identity $1/(1+d)-1+d= d^2-d^3/(1+d)$
and (\ref{o5}) yield
\[
\Biggl|\frac{1}{n}\sum_{j=1}^n
\biggl(\frac{r_j}{1+\zeta^{\top
}x_j}-r_j+r_j x_j^\top
\zeta \biggr) \Biggr| \leq \Biggl| \frac{1}{n}\sum_{j=1}^nr_j
\bigl(\zeta^{\top}x_j \bigr)^2 \Biggr| +
\frac{4}{3} \frac{1}{n}\sum_{j=1}^n|r_j|
\bigl|\zeta^{\top} x_j\bigr|^3
\]
for vectors $r_1,\dots,r_n$ of the same dimension.
Taking $r_j= S^{-1} x_j$, we derive with the help of (\ref{o2})
%
\[
\bigl|\zeta-S^{-1}\bar x\bigr|\leq \Biggl| \frac{1}{n}\sum
_{j=1}^nS^{-1} x_j \bigl(
\zeta^{\top}x_j \bigr)^2 \Biggr| + \frac{4}{3}
\frac{1}{n}\sum_{j=1}^n\bigl|S^{-1}x_j\bigr|\bigl|
\zeta^{\top}x_j\bigr|^3. 
\]
Using $|x|=\sup_{|v|=1} v^{\top}x$,
the Cauchy--Schwarz inequality,
(\ref{o4}) and (\ref{o6}) we bound the square of the first summand
of the right-hand side by
\[
\frac{1}{n}\sum_{j=1}^n \bigl(
\zeta^{\top} x_j \bigr)^4 \sup_{|v|=1}
v^{\top} S^{-1} v \le\frac{1}{\lambda} |\zeta|^4
x^{(4)}
\]
and the square of the second summand by
\[
\frac{16}{9\lambda^2} x_*^2 \zeta^{\top} S \zeta
\frac{1}{n}\sum_{j=1}^n \bigl(
\zeta^{\top
} x_j \bigr)^4 \le\frac{\Lambda}{9\lambda^2} |
\zeta|^4 x^{(4)}.
\]
Combining the above, we obtain
%
\begin{equation}
\label{o7} \bigl|\zeta-S^{-1}\bar x\bigr|^2 \leq2 \biggl(
\frac{1}{\lambda} + \frac
{\Lambda}{9\lambda^2} \biggr) |\zeta|^4
x^{(4)}.
\end{equation}
Using the inequality
$|2\log(1+t)-2t+t^2 -2t^3/3|\leq|t|^4/(2(1-|t|)^4)$ valid for $|t|<1$,
and then (\ref{o5}) we derive
\[
\Biggl| 2\sum_{j=1}^n \log \bigl(1+
\zeta^{\top}x_j \bigr) - 2n\zeta^{\top} \bar x + n
\zeta^{\top}S\zeta \Biggr| \leq\frac{2}{3} \Biggl| \sum
_{j=1}^n \bigl(\zeta^{\top}x_j
\bigr)^3 \Biggr| + \frac{1}{2} \biggl(\frac{4}{3}
\biggr)^4 \sum_{j=1}^n \bigl|
\zeta^{\top} x_j\bigr|^4. 
\]
With $\Delta= \zeta-S^{-1}\bar x$, we can write
$\zeta^{\top}S\zeta= \zeta^{\top}\bar x + \zeta^{\top}S \Delta$
and
$\zeta^{\top}\bar x = \bar x^{\top} S^{-1}\bar x + \Delta^{\top
}\bar x$,
and obtain the identity $2\zeta^{\top}\bar x - \zeta^{\top}S\zeta=
\bar x S^{-1}\bar x - \Delta^{\top}S\Delta$. Using this and (\ref{o7}),
we arrive at the bound
\[
\Biggl| 2\sum_{j=1}^n \log \bigl(1+
\zeta^{\top}x_j \bigr) - n\bar x^{\top
}S^{-1}
\bar x \Biggr| \leq n |\zeta|^3 x^{(3)} + n \biggl(
\frac{16}{9}+\frac{2\Lambda
}{\lambda} + \frac{2\Lambda^2}{9\lambda^2} \biggr) |
\zeta|^4 x^{(4)}. 
\]
In view of (\ref{o3}) and $\Lambda\ge\lambda$, this becomes
%
\begin{equation}
\label{o8a} \Biggl| 2\sum_{j=1}^n \log \bigl(1+
\zeta^{\top}x_j \bigr) - n\bar x^{\top
}S^{-1}
\bar x \Biggr| \leq\frac{n|\bar x|^3 x^{(3)}}{(\lambda-|\bar x|x_*)^3} +\frac{\Lambda^2}{\lambda^2}\frac{4n|\bar x|^4 x^{(4)}}{(\lambda
-|\bar x|x_*)^4}.
\end{equation}
If we bound $x^{(3)}$ by $x_*\Lambda$ and $x^{(4)}$ by $x_*^2 \Lambda
$ and use
(\ref{o4}), we obtain the bound
%
\begin{equation}
\label{o8} \Biggl| 2\sum_{j=1}^n \log \bigl(1+
\zeta^{\top}x_j \bigr) - n\bar x^{\top
}S^{-1}
\bar x \Biggr| \leq \biggl(\Lambda+ \frac{\Lambda^3}{\lambda^2} \biggr) \frac{n x_* |\bar x|^3}{(\lambda- |\bar x|x_*)^3}.
\end{equation}
%


Thus, we have proved the following result.

\begin{lemma}\label{lm.owen}
The inequality $\lambda> 5|\bar x| x_*$ implies that there is a unique
$\zeta$ in $\R^m$ satisfying
$1+\zeta^{\top}x_j>0$, $j=1,\dots,n$, and (\ref{o2}) to (\ref{o8}).
\end{lemma}

\section{Applications with random vectors} \label{app}

We shall now discuss implications of Lemma~\ref{lm.owen} to the case
when the
vectors $x_j$ are replaced by random vectors.
We are interested in the case when the dimension of the random vectors
increases with $n$.

Let $T_{n1},\dots,T_{nn}$ be $m_n$-dimensional random vectors.
With these random vectors we associate the empirical likelihood
\[
\el=\sup \Biggl\{\prod_{j=1}^n n
\pi_j\dvt\pi\in\prob, \sum_{j=1}^n
\pi_jT_{nj}=0 \Biggr\}.
\]
To study the asymptotic behavior of $\el$, we introduce
\[
T_n^*= \max_{1\leq j \leq n}|T_{nj}|,\qquad \bar
T_n = \frac{1}{n}\sum_{j=1}^nT_{nj},
\qquad T_n^{(\nu)}=\sup_{|u|=1} \Biggl|\frac{1}{n}\sum
_{j=1}^n \bigl(u^{\top}
T_{nj} \bigr)^{\nu} \Biggr|,\qquad
\nu=3,4
\]
and the matrix
\[
S_n = \frac{1}{n}\sum_{j=1}^nT_{nj}T_{nj}^{\top},
\]
and let $\lambda_n$ and $\Lambda_n$ denote the smallest
and largest eigenvalues of $S_n$,
\[
\lambda_n = \inf_{|u|=1} u^{\top}S_n u
\quad\mbox{and}\quad \Lambda_n = \sup_{|u|=1}
u^{\top}S_n u.
\]

We say a sequence $W_n$ of $m_n\times m_n$ dispersion matrices is \emph
{regular} if
the following condition holds,
\[
0 < \inf_n \inf_{|u|=1} u^{\top}W_n u
\leq\sup_n \sup_{|u|=1} u^{\top}W_n u <
\infty.
\]
We impose the following conditions.

\begin{enumerate}[(A4)]
\item[(A1)]
$m_n^{1/2}T_n^*= \mathrm{o}_p(n^{1/2})$.
\item[(A2)]
$n |\bar T_n|^2 = \mathrm{O}_p(m_n)$.
\item[(A3)]
There is a regular sequence of 
dispersion matrices $W_n$ such that
$
|S_n - W_{n}|_o = \mathrm{o}_p(m_n^{-1/2}).
$
\item[(A4)]
$m_n T_{n}^{(3)}= \mathrm{o}_p(n^{1/2})$ and $m_n^{3/2} T_n^{(4)}=\mathrm{o}_p(n)$.
\end{enumerate}

The first two conditions imply $T_n^*|\bar T_n| = \mathrm{o}_p(1)$,
the third condition implies that there are positive numbers $a<b$ such that
$P(a \leq\lambda_n \leq\Lambda_n \leq b) \to1$.
Thus, all three conditions imply that the probability of the event
$\{\lambda_n > 5 T_n^* |\bar T_n|\}$ tends to one.
Consequently, by Lemma~\ref{lm.owen}, there exists an $m_n$-dimensional
random vector $\hz$
which is uniquely determined on this event by the properties
$1+\hz^{\top} T_{nj}>0$, $j=1,\dots,n$,
and
%
\begin{equation}
\label{constraint} \frac{1}{n}\sum_{j=1}^n
\frac{T_{nj}}{1+\hz^{\top} T_{nj}}=0.
\end{equation}
On this event, we have $-2 \log\el= 2 \sum_{j=1}^n \log(1+\hz^{\top}T_{nj})$.
It follows from (A3) that $S_n$ is invertible except on an event
whose probability tends to zero.
It follows from (A2) and (A4) that
\[
n |\bar T_n|^3 T_n^{(3)} =
\mathrm{o}_p \bigl(m_n^{1/2} \bigr) \quad
\mbox{and} \quad n |\bar T_n|^4 T_n^{(4)}
= \mathrm{o}_p \bigl(m_n^{1/2} \bigr).
\]
%
Thus, under (A1)--(A4), the following expansion follows from (\ref{o8a})
%
\begin{equation}
\label{c2} -2 \log\el= n\bar T_n^{\top} S_n^{-1}
\bar T_n +\mathrm{o}_p \bigl(m_n^{1/2}
\bigr).
\end{equation}
%
From (A3), we can also derive the rate $|S_n^{-1}-W_n^{-1}|_o=\mathrm{o}_p(m_n^{-1/2})$.
Thus, if (A1)--(A4) hold, then (\ref{c2}) holds
with $S_n$ replaced by $W_{n}$,
%
\begin{equation}
\label{c2p} -2 \log\el= n\bar T_n^{\top}
W_{n}^{-1} \bar T_n +\mathrm{o}_p
\bigl(m_n^{1/2} \bigr).
\end{equation}
In view of the inequalities $T_n^{(3)} \le\Lambda_n T_n^*$ and
$T_n^{(4)} \le\Lambda_n (T_n^*)^2$, a sufficient condition for (A1)
and (A4) is given by
\renewcommand{\theequation}{B\arabic{equation}}
\setcounter{equation}{0}
\begin{equation}
\label{B1} m_n T_n^* = \mathrm{o}_p
\bigl(n^{1/2} \bigr).
\end{equation}
In view of the bound $(T_n^{(3)})^2 \le\Lambda_n T_n^{(4)}$, which is
a consequence of the
Cauchy--Schwarz
inequality, a sufficient condition for (A4) is given by
%
\begin{equation}
\label{B2} m_n^2 T_n^{(4)} =
\mathrm{o}_p(n).
\end{equation}
%

We first treat the case when the dimension $m_n$ does not increase with $n$.
In this case, (\ref{B1}) and (A2) are implied by $T_n^*=\mathrm{o}_p(n^{1/2})$
and $\bar T_n = \mathrm{O}_p(n^{-1/2})$, and (A3) is implied by the condition:
$S_n=W+\mathrm{o}_p(1)$ for some positive definite matrix $W$. Thus, we have the
following result.

\begin{theorem}\label{thm.1}
Let $m_n=m$ for all $n$. Suppose
\[
T_n^*=\mathrm{o}_p \bigl(n^{1/2} \bigr),\qquad
n^{1/2} \bar T_n \implies N(0,V) \quad\mbox{and}\quad
S_n = W+\mathrm{o}_p(1)
\]
for dispersion matrices $V$ and $W$, with $W$ positive definite.
Then $-2 \log\el$ converges in distribution
to $Z^{\top} V^{1/2}W^{-1} V^{1/2} Z$, where
the $m$-dimensional random vector $Z$ is standard normal.
For $V=W$, the limiting distribution is a chi-square distribution
with $m$ degrees of freedom.
\end{theorem}


If we replace 
$n^{1/2} \bar T_n \implies N(0,V)$ by
$n^{1/2} \bar T_n \implies U$ for some random variable $U$,
then the conclusion becomes $-2 \log\el$ converges in distribution
to $U^{\top} W^{-1} U$.
This version of the theorem yields Theorem 2.1 of Hjort \textit{et al.}
\cite{HMV09} without their (A0).

Theorem~\ref{thm.1} does not require the independence of the random vectors
$T_{n,1},\dots,T_{n,n}$.
This is important when dealing with estimated constraint functions as
we shall see below.

Suppose the condition in the theorem hold with $V=W$.
Under a contiguous alternative, one typically has
$n^{1/2}\bar T_n \implies N(\mu,V)$ for some $\mu$ different from zero,
but retains the other conditions. In this case, $-2 \log\el$ has a
limiting chi-square distribution with
$m$ degrees of freedom and non-centrality parameter $|V^{-1/2}\mu|$.

Let us address some applications of Theorem~\ref{thm.1}.
For this discussion, we let $Z_1,\dots,Z_n$ be independent copies of a
$k$-dimensional
random vector $Z$ with distribution $Q$
and let $w$ be a measurable function from $\R^k$ into $\R^m$
such that $E[w(Z)]=\int w \, \mathrm{d}Q=0$ and $W=E[w(Z)w^{\top}(Z)]=\int
ww^{\top}\, \mathrm{d}Q$
is positive definite. Let us first look at the empirical likelihood
\[
\mathscr{R}_{n1} =\sup \Biggl\{\prod_{j=1}^n
n\pi_j\dvt\pi\in\prob, \sum_{j=1}^n
\pi_jw(Z_j)=0 \Biggr\}.
\]
It follows from Owen that $-2\log\mathscr{R}_{n1}$ has a limiting
chi-square distribution with
$m$ degrees of freedom. This also follows from Theorem~\ref{thm.1}
applied with
$T_{nj}=w(Z_j)$.
Indeed, the first condition follows from the inequality
\renewcommand{\theequation}{\arabic{section}.\arabic{equation}}
\setcounter{equation}{3}
\begin{equation}
\label{fm1} 
P \Bigl(
\max_{1\leq j \leq n}\bigl|w(Z_j)\bigr| > \epsilon n^{1/2} \Bigr) \le
\frac{1}{\epsilon^2} E\bigl[\bigl|w(Z)\bigr|^2 \1 \bigl[\bigl|w(Z)\bigr|> \epsilon
n^{1/2} \bigr]\bigr]
\end{equation}
and the Lebesgue dominated convergence theorem; the central limit
theorem yields
the second condition with $V=W$;
the third condition
%
\begin{equation}
\label{fm3} \frac{1}{n}\sum_{j=1}^nw(Z_j)w^{\top}(Z_j)
= W +\mathrm{o}_p(1)
\end{equation}
follows from the weak law of large numbers.
This shows that Owen's result is a special case of our result.

Now consider the empirical likelihood
\[
\hat{\mathscr{R}}_{n1} = \sup \Biggl\{\prod
_{j=1}^n n\pi_j\dvt\pi\in \prob, \sum
_{j=1}^n \pi_j\hat
w(Z_j)=0 \Biggr\},
\]
where $\hat w$ is an estimator of $w$ based on the observations
$Z_1,\dots,Z_n$
which is consistent in the following sense,
%
\begin{equation}
\label{fm4} \frac{1}{n}\sum_{j=1}^n\bigl|
\hat w(Z_j)-w(Z_j)\bigr|^2 =
\mathrm{o}_p(1).
\end{equation}
Then $-2\log\hat{\mathscr{R}}_{n1}$ has a limiting chi-square distribution
with $m$ degrees of freedom if also
%
\begin{equation}
\label{fm5} n^{-1/2}\sum_{j=1}^n
\hat w(Z_j)= n^{-1/2}\sum_{j=1}^nw(Z_j)+
\mathrm{o}_p(1)
\end{equation}
holds.
To see this, we verify the assumptions of Theorem~\ref{thm.1} with
$T_{nj}=\hat w(Z_j)$.
The first condition follows from (\ref{fm1}), (\ref{fm4}) and the inequality
\[
T_n^* \le\max_{1\leq j \leq n}\bigl|w(Z_j)\bigr| + \Biggl( \sum
_{j=1}^n \bigl|\hat w(Z_j)-w(Z_j)\bigr|^2
\Biggr)^{1/2}.
\]
The central limit theorem, Slutsky's theorem and (\ref{fm5}) yield
the second condition with $V=W$.
The third condition follows from (\ref{fm3}), (\ref{fm4})
and the inequality (\ref{onorm.diff}).

The requirement (\ref{fm5}) is rather strong. One often only derives
%
\begin{equation}
\label{fm6} n^{-1/2}\sum_{j=1}^n
\hat w(Z_j)= n^{-1/2}\sum_{j=1}^nv(Z_j)
+ \mathrm{o}_p(1)
\end{equation}
for some function $v$ satisfying $E[v(Z)]=0$ and $E[|v(Z)|^2] <\infty$.
Under (\ref{fm4}) and (\ref{fm6}), $-2\log\hat{\mathscr{R}}_{n1}$
has limiting
distribution as given in Theorem~\ref{thm.1} with $V$ the dispersion
matrix of $v(Z)$.
This follows from Theorem~\ref{thm.1} whose assumptions are now
verified as above.

In situations when $w(Z)=u(Z,\eta)$ for some $q$-dimensional nuisance parameter
$\eta$ and $\hat w(Z)= u(Z,\hat\eta)$ for some estimator $\hat\eta
$ of $\eta$,
one typically has $v(Z)=w(Z)+D\psi(Z)$, where the $m\times q$ matrix
$D$ is the
derivative of the map $t \mapsto E[u(Z,\eta+t)]$ at $t=0$, and $\psi$
is the
influence function of $\hat\eta$.

We now address the case when $m_n$ increases with the sample size.

\begin{theorem}\label{thm.2}
Let \emph{(A1)--(A4)} hold.
Suppose that $m_n$ increases with $n$ to infinity and that
there are $m_n\x m_n$ dispersion matrices $V_n$
such that $m_n/\trace(V_n^2)=\mathrm{O}(1)$ and
%
\begin{equation}
\label{twt} 
\bigl(n \bar T_n^{\top}
W_{n}^{-1}\bar T_n - \trace(V_n)
\bigr) /\sqrt{2 \trace \bigl(V_n^2 \bigr)}
\implies N(0,1).
\end{equation}
Then we have
%
\begin{equation}
\label{lR} 
\bigl(-2 \log\el-\trace(V_n)
\bigr) /\sqrt{2 \trace \bigl(V_n^2 \bigr)}
\implies N(0,1).
\end{equation}
\end{theorem}

\begin{pf} We have already seen that (A1)--(A4) imply (\ref{c2p}).
It follows from (\ref{c2p}) and $m_n/\trace(V_n^2)=\mathrm{O}(1)$ that
the difference of the left-hand sides of\vadjust{\goodbreak} (\ref{twt}) and (\ref{lR}) converge
to zero in probability. Thus, the desired (\ref{lR}) follows from
(\ref{twt}) and
Slutsky's theorem.
\end{pf}

Of special interest is the case when $V_n$ is the $m_n\x m_n$ identity matrix
$I_{m_n}$.
Then $\trace(V_n)=\trace(V_n^2)=m_n$ and (\ref{lR}) simplifies to
(\ref{hmv}).
Sufficient conditions for (\ref{twt}) are given by Peng and Schick~\cite{PS12}.


\section{Main results} \label{main}
In this section, we assume that $(\Z,\mathscr{S})$ is a measurable space,
that $Z_1,\dots,Z_n$ are independent copies of the
$\Z$-valued random variable $Z$ with distribution $Q$,
and that $m_n$ is a positive integer that tends to infinity with $n$.
We let $w_n$ denote a measurable function from $\Z$ to $\R^{m_n}$
such that
$\int w_n\, \mathrm{d}Q=0$ and $\int|w_n|^2 \,\mathrm{d}Q$ is finite.

We first study
\[
\el= \sup \Biggl\{\prod_{j=1}^n n
\pi_j\dvt\pi\in\prob, \sum_{j=1}^n
\pi_jw_n(Z_j)=0 \Biggr\}.
\]
Our goal is to show (\ref{hmv}).
To this end, we set
\[
\bar w_n= \frac{1}{n}\sum_{j=1}^nw_n(Z_j),
\qquad \bar W_n = \frac{1}{n}\sum_{j=1}^nw_n(Z_j)w_n^{\top}(Z_j),
\qquad W_n = \int w_nw_n^{\top} \,
\mathrm{d}Q
\]
and introduce the following condition.

\begin{enumerate}
\item[(C1)]
The sequence $W_n$ is regular.
\end{enumerate}

Motivated by the results in Peng and Schick~\cite{PS12}, we call a sequence
$v_n$ of measurable functions from $\Z$ to $\R$ \emph{Lindeberg} if
%
\begin{equation}
\label{L} \int|v_n|^2 \1\bigl[|v_n|> \epsilon
\sqrt{n}\bigr]\, \mathrm{d}Q \to0,\qquad \epsilon>0.
\end{equation}
The following are easy to check.
If the sequences $u_n$ and $v_n$ are Lindeberg, so are the sequences
$\max\{|u_n|,|v_n|\}$ and $u_n +v_n$. If the sequence $v_n$ is
Lindeberg and $|u_n|\le|v_n|$,
then the sequence $u_n$ is also Lindeberg. We also need the following
properties.
\begin{enumerate}
\item[(L1)]
If $v_n$ is Lindeberg, then one has the rate
$
\max_{1\leq j \leq n}|v_n(Z_j)| = \mathrm{o}_p(n^{1/2}).
$
\item[(L2)]
If $\int|v_n|^r \,\mathrm{d}Q=\mathrm{o}(n^{r/2-1})$ for some $r>2$, then $v_n$ is Lindeberg.
\end{enumerate}
The first statement follows from an inequality similar to (\ref{fm1}),
the second from Remark~1 in Peng and Schick~\cite{PS12}.
%

To show (\ref{hmv}), we apply Theorem~\ref{thm.2} with $T_{nj}=w_n(Z_j)$.
In the presence of (C1), the conditions
(\ref{twt}) and (A1)--(A4) of this theorem are implied by
{\makeatletter
\setvaluelist{myfoo}{k,D0,D1,D2,D3,D4}
\def\theequation{\getitemvalue{myfoo}{\the\c@equation}}
\makeatother
\begin{eqnarray}
\label{D0}\bigl(n \bar w_n^{\top}W_n^{-1}
\bar w_n-m_n \bigr) / \sqrt{2m_n} &
\implies& N(0,1),
\\
\label{D1}\max_{1\leq j \leq n}m_n^{1/2}
\bigl|w_n(Z_j)\bigr| &=& \mathrm{o}_p
\bigl(n^{1/2} \bigr),
\\
\label{D2}n|\bar w_n|^2 &=& \mathrm{O}_p(m_n),
\\
\label{D3}|\bar W_n-W_n|_o &=&
\mathrm{o}_p \bigl(m_n^{-1/2} \bigr),
\\
\label{D4}\sup_{|u|=1} \frac{m_n^2}{n} \sum
_{j=1}^n \bigl|u^{\top} w_n(Z_j)\bigr|^4
&=& \mathrm{o}_p(n).
\end{eqnarray}}%
%
By part (c) of Corollary 3 in Peng and Schick~\cite{PS12},
(\ref{D0}) follows if the function $|W_n^{-1/2}w_n|$ is Lindeberg.
In the presence of (C1), the latter condition is equivalent to
$|w_n|$ being Lindeberg.
By (L1), a sufficient condition for (\ref{D1}) is that $m_n^{1/2}|w_n|$ is Lindeberg.
It follows from (C1) that $\trace(W_n)\le B m_n$ for some constant
$B$. Thus (C1) implies
$E[n|\bar w_n|^2] = \trace(W_n)=\mathrm{O}(m_n)$ and hence (\ref{D2}).
In view of (C1), a sufficient condition for (\ref{D3}) is that $m_n|w_n|$ is
Lindeberg.
To see this, fix $\epsilon>0$ and let $\bar W_{n,1}$ and $\bar
W_{n,2}$ be the matrices
obtained by replacing in the definition of $\bar W_{n}$
the function $w_n$ by $v_n=w_n \1[|m_n w_n| \leq\epsilon\sqrt{n}]$ and
$w_n-v_n=w_n \1[|m_n w_n|> \epsilon\sqrt{n}]$, respectively.
Then we find
\[
n E \bigl[\bigl|\bar W_{n,1}-E[\bar W_{n,1}]\bigr|^2 \bigr]
\leq E \bigl[|v_n|^4(Z) \bigr]
\leq \frac{\epsilon^2 n}{m_n^2} E \bigl[|w_n|^2(Z) \bigr] \leq
\frac{\epsilon^2 n B m_n}{m_n^{2}},
\]
\[
P(\bar W_{n,2}\neq0)\leq P \Bigl(\max_{1\leq j \leq n}\bigl|m_n
w_n(Z_j)\bigr|>\epsilon\sqrt{n} \Bigr)\to0
\]
and using (\ref{onorm.cov})
\[
\bigl|E[\bar W_{n,2}]\bigr|_o \leq E \bigl[|w_n|^2(Z)
\1 \bigl[\bigl|m_n w_n(Z)\bigr|>\epsilon\sqrt{n} \bigr] \bigr] =
\mathrm{o} \bigl(m_n^{-2} \bigr).
\]
The above inequalities show that (C1) and $m_n|w_n|$ is Lindeberg imply
statement (\ref{D3}).
The latter condition also implies (\ref{B1}) and hence (\ref{D1}) and (\ref{D4}), the latter
in the presence of (C1).
Thus, we have the following result.

\begin{theorem}\label{thm.3}
Suppose \emph{(C1)} holds and the sequence $m_n|w_n|$ is Lindeberg.
Then (\ref{hmv}) holds as $m_n$ tends to infinity with $n$.
\end{theorem}

From this, simple calculations and the property (L2) we immediately
derive the following corollaries.

\begin{corollary}\label{cor.3a}
Suppose \emph{(C1)} holds and $|w_n|\le\sqrt{m_n} B$ for some constant $B$.
Then (\ref{hmv}) holds if $m_n^3=\mathrm{o}(n)$.
\end{corollary}

\begin{corollary}\label{cor.3b}
Suppose \emph{(C1)} holds and $\int|w_n|^r \,\mathrm{d}Q= \mathrm{O}(m_n^{r/2})$ for some $r>2$.
Then (\ref{hmv}) holds if $m_n^{3r/(r-2)} =\mathrm{o}(n)$.
\end{corollary}


These two corollaries give the conclusions in Theorem 4.1 in Hjort \textit{et al.}~\cite{HMV09}
under slightly weaker conditions in the case of Corollary~\ref{cor.3b}.
We now present some additional\vadjust{\goodbreak} results that allow for larger $m_n$
if $r$ is small.
For example, if $r=4$,
Corollary~\ref{cor.3b} requires $m_n^6=\mathrm{o}(n)$, while Theorem \ref
{thm.3b} below
allows $m_n^4=\mathrm{o}(n)$.
For $r=3$,
Corollary~\ref{cor.3b} requires $m_n^9=\mathrm{o}(n)$, while Theorem \ref
{thm.3d} below
allows $m_n^6=\mathrm{o}(n)$.

\begin{theorem}\label{thm.3b}
Suppose \emph{(C1)} holds and $\int|w_n|^4\, \mathrm{d}Q= \mathrm{O}(m_n^{2})$.
Then (\ref{hmv}) holds if \mbox{$m_n^4=\mathrm{o}(n)$}.
\end{theorem}

\begin{pf}
Using (L2) and $m_n^4=\mathrm{o}(n)$ we derive that $m_n^{1/2} |w_n|$ is Lindeberg.
This latter condition and (C1) imply (\ref{D1})--(\ref{D2}) as shown prior to
Theorem~\ref{thm.3}.
Next we calculate $nE[|\bar W_n -W_n|^2] \leq E[|w_n|^4(Z)]= \mathrm{O}(m_n^2)$.
This yields (\ref{D3}) in view of
$|\bar W_n -W_n|_o \le\break|\bar W_n-W_n|= \mathrm{O}_p(m_n/\sqrt{n})$ and $m_n^4=\mathrm{o}(n)$.
Finally, we have (\ref{D4}) as the left-hand side of (\ref{D4}) is bounded by
\[
\frac{m_n^2}{n} \sum_{j=1}^n
\bigl|w(Z_j)\bigr|^4 = \mathrm{O}_p
\bigl(m_n^4 \bigr) =\mathrm{o}_p(n).
\]
Thus, (\ref{D0})--(\ref{D4}) hold and we obtain the desired result from Theorem
\ref{thm.2}.
\end{pf}

\begin{theorem}\label{thm.3d}
Suppose \emph{(C1)} holds and $\int|w_n|^r\, \mathrm{d}Q= \mathrm{O}(m_n^{r/2})$ for some $2< r < 4$.
Then (\ref{hmv}) holds if $m_n^{2r/(r-2)}=\mathrm{o}(n)$.
\end{theorem}

\begin{pf}
There is a constant $B$ such that $\int|w_n|^r \, \mathrm{d}Q \le B m_n^{r/2}$.
In view of (L2) and the properties of $m_n$, we derive that $m_n^{1/2}
|w_n|$ is Lindeberg.
This condition and (C1) imply (\ref{D0})--(\ref{D2}). It follows from (\ref{D1}), the
moment condition
on $w_n$, and the properties of $m_n$ that
\begin{eqnarray*}
\frac{m_n^2}{n} \sum_{j=1}^n
\bigl|w(Z_j)\bigr|^4 &\le& \frac{m_n^2}{n} \sum
_{j=1}^n \bigl|w(Z_j)\bigr|^r
\max_{1\leq j \leq
n}\bigl|w_n(Z_j)\bigr|^{4-r}
\\
&=& \mathrm{o}_p \bigl( m_n^2
m_n^{r/2} (n/m_n)^{(4-r)/2} \bigr) =
\mathrm{o}_p \bigl(m_n^r n^{(4-r)/2}
\bigr)=\mathrm{o}_p(n).
\end{eqnarray*}
This establishes (\ref{D4}). Finally, (\ref{D3}) follows as we have
$|\bar W_n - W_n|_o = \mathrm{o}_p(m_n^{-1})$.
To prove the latter, we mimic the argument prior to Theorem~\ref{thm.3}
used to verify (\ref{D3}) if $m_n|w_n|$ is Lindeberg.
But now $|w_n|\1[m_n^{1/2}|w_n| \le\sqrt{n}]$
plays the role of $v_n$.
For the corresponding matrices $\bar W_{n1}$ and $\bar W_{n2}$, we have
\begin{eqnarray*}
m_n^2 E\bigl[\bigl|\bar W_{n1}-E [\bar
W_{n1}]\bigr|^2 \bigr] &\le&\frac{m_n^2}{n} \biggl(
\frac{n}{m_n} \biggr)^{(4-r)/2} B m_n^{r/2} \le
\frac{B m_n^{r}}{n^{r/2-1}} \to0,
\\
P(\bar W_{n2}\ne0)&\le& P \Bigl(\max_{1\leq j \leq n}m_n^{1/2}
\bigl|w_n(Z_j)\bigr| > n^{1/2} \Bigr) \to0,
\\
m_n \bigl|E[\bar W_{n,2}]\bigr|_o &\le&\int
\frac{m_n^{r/2}|w_n|^r}{n^{(r-2)/2}}\, \mathrm{d}Q \le\frac{B m_n^{r}}{n^{r/2-1}} \to0.
\end{eqnarray*}
Consequently, (\ref{D0})--(\ref{D4}) hold and the desired result follows.\vadjust{\goodbreak}
\end{pf}


Now, we study
\[
\hel= \sup \Biggl\{\prod_{j=1}^n n
\pi_j\dvt\pi\in\prob, \frac
{1}{n}\sum
_{j=1}^n\pi_j \hat w_n(Z_j)=0
\Biggr\},
\]
where $\hat w_n$ is an estimator of $w_n$. Let us set
\[
\hat W_n = \frac{1}{n}\sum_{j=1}^n
\hat w_n(Z_j)\hat w_n^{\top}(Z_j).
\]

\begin{theorem}\label{thm.hw}
Suppose \emph{(C1)} holds and assume we have the expansions
\setcounter{equation}{1}
\begin{eqnarray}
\label{a1} m_n \max_{1\leq j \leq n}\bigl|\hat w_n(Z_j)\bigr|
&=& \mathrm{o}_p \bigl(n^{1/2} \bigr),
\\
\label{a2} |\hat W_n - W_n|_o &=&
\mathrm{o}_p \bigl(m_n^{-1/2} \bigr),
\\
\label{a3} \frac{1}{n}\sum_{j=1}^n
\hat w_n(Z_j) &=& \frac{1}{n}\sum
_{j=1}^nv_n(Z_j) +
\mathrm{o}_p \bigl(n^{-1/2} \bigr)
\end{eqnarray}
for some measurable function $v_n$ from $S$ into $\R^{m_n}$ such that
$\int v_n \,\mathrm{d}Q=0$ and $|v_n|$ is Lindeberg.
Furthermore, assume that the dispersion matrix
\[
U_n=W_n^{-1/2} \int v_nv_n^{\top}
\, \mathrm{d}Q W_n^{-1/2}
\]
of $W_n^{-1/2} v_n(Z)$ satisfies $|U_n|_o=\mathrm{O}(1)$
and $m_n/\trace(U_n^2)$ is bounded.
Then, as $m_n$ tends to infinity with $n$,
$
(-2 \log\hel- \trace(U_n))/
{\sqrt{2 \trace(U_n^2)}}
$
is asymptotically standard normal.
\end{theorem}

\begin{pf}
Set $\xi_{nj}= W_n^{-1/2} v_n(Z_j)$, and introduce the averages
$\bar v_n= \frac{1}{n}\sum_{j=1}^nv_n(Z_j)$ and $\bar T_n= \frac
{1}{n}\sum_{j=1}^n\hat w_n(Z_j)$.
It follows from (C1) that $|W_n^{1/2}|_o+|W_n^{-1/2}|_o= \mathrm{O}(1)$.
Using this and the Lindeberg property of $|v_n|$, we derive
%
\begin{equation}
\label{li.0} L_n(\epsilon)=E \bigl[|\xi_{n,1}|^2
\1\bigl[|\xi_{n,1}|> \epsilon\sqrt{n}\bigr] \bigr] \to0,\qquad \epsilon>0.
\end{equation}
We have
$\trace(U_n)/\trace(U_n^2) \leq|U_n|_o m_n/\trace(U_n^2)=\mathrm{O}(1)$.
From $m_n/\trace(U_n^2) =\mathrm{O}(1)$ we conclude
$\trace(U_n^2) \to\infty$.
Thus, Theorem 2 in Peng and Schick~\cite{PS12} yields that
$
(n \bar v_n W_n^{-1} \bar v_n - \trace(U_n) ) /
{\sqrt{2\trace(U^2_n)} }
$
is asymptotically standard normal.
From this, (C1), $\trace(U_n)=\mathrm{O}(m_n)$ and $\trace(U_n^2)\le|U_n|_0^2m_n$
we conclude
$n|\bar v_n|^2 = \mathrm{O}_p(m_n)$.
With the help of (\ref{a3}) and the assumption $m_n/\trace(U^2_n)=\mathrm{O}(1)$,
we then derive
$n|\bar T_n|^2 = O_p(m_n)$ and that
$
(n \bar T_n W_n^{-1} \bar T_n - \trace(U_n) ) /
{\sqrt{2\trace(U^2_n)}} 
$
is asymptotically standard normal.
Thus in view of (\ref{B1}), conditions (A1)--(A4) hold with $T_{nj}=\hat w_n(Z_j)$,
and the desired result follows from Theorem~\ref{thm.2}.
\end{pf}


Let us first mention the special case when $v_n=w_n$.
In this case, $U_n$ equals $I_{m_n}$ and $\trace(U_n)=\trace(U_n^2)=m_n$.

\begin{corollary}\label{cor.1}
Suppose \emph{(C1)}, (\ref{a1}) and (\ref{a2}) hold,
$|w_n|$ is Lindeberg, and the following expansion is valid,
%
\begin{equation}
\label{a3w1} \frac{1}{n}\sum_{j=1}^n
\hat w_n(Z_j) = \frac{1}{n}\sum
_{j=1}^nw_n(Z_j) +
\mathrm{o}_p \bigl(n^{-1/2} \bigr).
\end{equation}
Then $(-2 \log\hel- m_n)/{\sqrt{2m_n}}$
is asymptotically standard normal.
\end{corollary}

Next, we treat $v_n=w_n-A_n\psi$ with $A_n$ and $\psi$ as
in the next condition.

\begin{enumerate}[(C2)]
\item[(C2)]
There is a measurable function $\psi$ from $\Z$ into $\R^q$ satisfying
$\int\psi\, \mathrm{d}Q=0$ and $\int\psi\psi^{\top}\, \mathrm{d}Q=I_q$ such that,
with $A_n =\int w_n \psi^{\top}\, \mathrm{d}Q$, the expansion,
\[
\frac{1}{n}\sum_{j=1}^n
\hat w_n(Z_j) = \frac{1}{n}\sum
_{j=1}^nw_n(Z_j)-
A_n \psi(Z_j) +\mathrm{o}_p
\bigl(n^{-1/2} \bigr),
\]
and 
the convergence, $ \trace(A_n^{\top}W_n^{-1} A_n) \to q$, hold.
\end{enumerate}

\begin{corollary}\label{cor.2}
Suppose \emph{(C1)}, \emph{(C2)}, (\ref{a1}) and (\ref{a2}) hold,
and $|w_n|$ is Lindeberg.
Then $(-2 \log\hel- m_n+q)/{\sqrt{2(m_n-q)}}$
is asymptotically standard normal.
\end{corollary}

\begin{remark}\label{rem.fc}
Suppose that $w_n$ is the vector formed by the first $m_n$ elements
of an orthonormal basis $u_1,u_2,\dots$ for $L_{2,0}(Q)$.
Then the $\nu$th column of the matrix $A_n$ is formed by the
first $m_n$ Fourier coefficients of the $\nu$th component of $\psi$
with respect to this basis. In this case, we have the identity
\[
\trace \bigl(A_n^{\top} W_n^{-1}A_n
\bigr) = \trace \bigl(A_n^{\top}A_n \bigr) = \sum
_{\nu=1}^q \sum_{k=1}^{m_n}
\biggl(\int\psi_{\nu} u_k\, \mathrm{d}Q
\biggr)^2
\]
and obtain under the assumptions $\int\psi\, \mathrm{d}Q=0$ and
$\int\psi\psi^{\top}\, \mathrm{d}Q=I_q$ the convergence
\[
\trace\bigl(A_n^{\top} W_n^{-1}A_n\bigr) \to\int|\psi|^2\, \mathrm{d}Q =q.
\]
\end{remark}

In our goodness-of-fit examples, 
the following condition holds.

\begin{enumerate}
\item[(C3)]
There is a constant $B$ such that $|w_n|\le B \sqrt{m_n}$ and $|\hat
w_n| \le B \sqrt{m_n}$.
\end{enumerate}

Under this condition, the rate $m_n^3/n \to0$ implies (\ref{a1}),
the Lindeberg property of $m_n |w_n|$, and (\ref{D3}).
Sufficient conditions for (\ref{a2}) can now be given directly or by verifying
%
\begin{equation}
\label{sa2a} |\hat W_n -\bar W_n|_o =
\mathrm{o}_p \bigl(m_n^{-1/2} \bigr).
\end{equation}
In view of the inequality (\ref{onorm.diff}), a sufficient condition for
the latter is
%
\begin{equation}
\label{sa2} D_n = \frac{1}{n}\sum
_{j=1}^n\bigl| \hat w_n(Z_j)-w_n(Z_j)\bigr|^2
= \mathrm{o}_p \bigl(m_n^{-1} \bigr).
\end{equation}

Thus, we have the following results.

\begin{corollary}\label{cor.3}
Suppose \emph{(C1)}, \emph{(C3)}, $m_n^3=\mathrm{o}(n)$,
and one of (\ref{a2}), (\ref{sa2a}), (\ref{sa2}) hold.
Then
\emph{(i)} (\ref{a3w1}) implies that $(-2 \log\hel- m_n)/{\sqrt{2m_n}}$
is asymptotically standard normal, while
%
\emph{(ii)} \emph{(C2)} implies that $(-2 \log\hel- m_n+q)/{\sqrt{2(m_n-q)}}$
is asymptotically standard normal.
\end{corollary}

\begin{remark}\label{rem.c0}
The conditions in Theorem~\ref{thm.hw} are based on the
sufficient condition (\ref{B1}) for (A1) and (A4).
Working with (A1) and (\ref{B2}) instead,
we see that (\ref{a1}) can be replaced
by the conditions,
\[
m_n^{1/2} \max_{1\leq j \leq n}\bigl|\hat
w_n(Z_j)\bigr| = \mathrm{o}_p
\bigl(n^{1/2} \bigr) 
\quad\mbox{and}
\quad\frac{m_n^2}{n} \sum_{j=1}^n \bigl|\hat
w_n(Z_j)\bigr|^4 = \mathrm{o}_p(n).
\]
With $D_n$ as in (\ref{sa2}), we derive the bounds
\begin{eqnarray*}
\max_{1\leq j \leq n}\bigl|\hat w_n(Z_j)\bigr| &\leq&
\max_{1\leq j \leq
n}\bigl|w_n(Z_j)\bigr| + (nD_n)^{1/2},
\\
\frac{m_n^2}{n} \sum_{j=1}^n \bigl|\hat
w_n(Z_j)\bigr|^4 &\le& \frac{8m_n^2}{n} \sum
_{j=1}^n \bigl|w_n(Z_j)\bigr|^4
+ 8m_n^2 nD_n^2.
\end{eqnarray*}
Here we used that $(a+b)^4 \le8(a^4+b^4)$ for nonnegative $a$ and $b$.
Assume now that $\int|w_n|^4 \,\mathrm{d}Q = \mathrm{O}(m_n^2)$ and that $m_n^4/n \to0$.
Then we have (\ref{D1}) and (\ref{D3}) as shown in the proof of Theorem~\ref{thm.3b}
and obtain the above two conditions
and (\ref{a2}) from (\ref{sa2}).
\end{remark}

\begin{corollary}\label{cor.5}
Suppose \emph{(C1)},
(\ref{sa2}), $\int|w_n|^4 \, \mathrm{d}Q=\mathrm{O}(m_n^2)$
and $m_n^4=\mathrm{o}(n)$ hold. Then
\emph{(i)} (\ref{a3w1}) implies that
$(-2 \log\hel- m_n)/{\sqrt{2m_n}}$
is asymptotically standard normal, while
\emph{(ii)} \emph{(C2)} implies that
$(-2 \log\hel- m_n+q)/{\sqrt{2(m_n-q)}}$
is asymptotically standard normal.\looseness=1
\end{corollary}

\begin{remark}\label{rem.alt}
Let us now describe the behavior of $-2\log\hel$ under a local alternative.
For this, we follow Remarks 6 and 7 in Peng and Schick~\cite{PS12}.
As there let $h$ be a measurable function satisfying
$\int h \,\mathrm{d}Q=0$ and $\int h^2 \,\mathrm{d}Q<\infty$ and let
$Q_{n,h}$ be a distribution satisfying\looseness=1
%
\begin{equation}
\label{hd1} \int\bigl|n^{1/2}(\sqrt{\mathrm{d}Q_{n,h}}-\sqrt{
\mathrm{d}Q}) -(1/2) h \sqrt{\mathrm{d}Q}\bigr|^2 \to0.
\end{equation}
Then the product measures $Q_{n,h}^n$ and $Q^n$ are mutually contiguous.
All results in this section obtain the expansion
%
\begin{equation}
\label{l1} -2 \log\hel- \Biggl|n^{-1/2}\sum_{j=1}^nu_n(Z_j)\Biggr|^2
= \mathrm{o}_p \bigl( m_n^{1/2} \bigr)
\end{equation}
for some measurable function $u_n$ from $\Z$ into $\R^{m_n}$
with the properties
\mbox{$\int u_n \,\mathrm{d}Q=0$}, $\int|u_n|^2\, \mathrm{d}Q = \mathrm{O}(m_n)$, $|u_n|$ is Lindeberg,
and the matrix $U_n= \int u_nu_n^{\top} \,\mathrm{d}Q$ satisfies
$|U_n|_o=\mathrm{O}(1)$ and $m_n/\trace(U_n^2)=\mathrm{O}(1)$. For example, in Theorem
\ref{thm.hw}
one has $u_n = W_n^{-1/2}v_n$.
By contiguity, one has the expansion (\ref{l1})
even if $Z_1,\dots,Z_n$ are independent with distribution $Q_{n,h}$.
Under this distributional assumption, one has
\[
\Biggl(\Biggl|n^{-1/2}\sum_{j=1}^nu_n(Z_j)\Biggr|^2
- \bigl|\mu_n(h)\bigr|^2 - \trace (U_n) \Biggr)\big / {
\sqrt{2\trace \bigl(U_n^2 \bigr)}} \implies
N(0,1)
\]
with $\mu_n(h)= \int u_n h \,\mathrm{d}Q$. Thus, under
the local alternative $Q_{n,h}$ one has
\[
\bigl(-2\log\hel-\bigl|\mu_n(h)\bigr|^2 - \trace(U_n)
\bigr) / {\sqrt{2\trace \bigl(U_n^2 \bigr)}}
\implies N(0,1).
\]
If $U_n=I_{m_n}$, this simplifies to
$(-2 \log\hel- |\mu_n(h)|^2 - m_n)/\sqrt{2m_n}\implies N(0,1)$
and may be interpreted as $-2\log\el$ being approximately
a non-central chi-square random variable with $m_n$ degrees of freedom
and non-centrality parameter $|\mu_n(h)|$.
\end{remark}

\section{Details for the examples} \label{det}

In this section, we use the results of the previous section to provide
the details for the examples of Sections~\ref{ex}.
In all examples, the components of $w_n$ are orthonormal and
uniformly bounded, so that (C1) and (C3) hold with $W_n=I_{m_n}$.
We begin with a technical lemma.

\begin{lemma}\label{lm.tech}
Let $(S_1,T_1),\dots,(S_n,T_n)$ be independent copies of the bivariate
random vector $(S,T)$,
where $T$ has a continuous distribution function $H$ and $E[S|T]=0$ and
$\sigma^2(T)=E[S^2|T]$ is bounded (by say $B$) and bounded away from
zero (by say $b$),
Let $\H$ denote the empirical distribution function based on
$T_1,\dots,T_n$.
Set $u_r=(1,\phi_1,\dots,\phi_r)^{\top}$, $D_j= u_r(\H(T_j))-u_r(H(T_j))$,
and $M=E[S^2 u_r(H(T)) u_r^{\top}(H(T)]$.
Then we have the following inequalities
%
\begin{equation}
\label{tech0} b \le v^{\top} M v \le B,\qquad v\in\R^{1+r},
|v|=1,
\end{equation}
%
\begin{eqnarray}
\nonumber
\\[-25pt]
\label{tech1} \Biggl|\frac{1}{n}\sum_{j=1}^nu_r
\bigl(\H(T_j) \bigr)u_r^{\top} \bigl(\H
(T_j) \bigr)-I_{1+r} \Biggr|^2 &\le&
\frac{16\uppi^2 r^2(1+r)^2}{n^2} \qquad\mbox{a.s.},
\\
\label{tech2} \frac{1}{n}\sum_{j=1}^nE
\bigl[|S_j D_j|^2 \bigr] &\le&
\frac{1}{n} \sum_{j=1}^nB E
\bigl[|D_j|^2 \bigr] \le \frac{B\uppi^2r^3}{n},
\\
\label{tech3} E \Biggl[ \Biggl|n^{-1/2}\sum_{j=1}^nS_j
D_j \Biggr|^2 \Biggr] &\le&\frac
{2B\uppi^2r^3}{n}. 
\end{eqnarray}
Moreover, if $E[S^4]$ is finite, then we have the bound
\[
E\bigl[S^4 \bigl|u_r \bigl(H(T) \bigr)\bigr|^4\bigr)
\le(1+2r)^2 E \bigl[S^4 \bigr].
\]
\end{lemma}

\begin{pf}
The last inequality follows from the bound $|u_r|^2 \le1+2r$.
The inequality (\ref{tech0}) is an easy consequence of $b \le\sigma^2(T)\le B$.
Conditioning on $T_1,\dots,T_n$ shows that
the left-hand side of (\ref{tech3}) is bounded by the left-hand side of
(\ref{tech2})
and yields the first inequality in (\ref{tech2}). Since $|\phi_k'|\le
\sqrt{2}\uppi k$,
we obtain $|D_j|^2 \le2\uppi r^3 (\H(T_j)-H(T_j))^2$. It is easy to
check that
$E[(\H(T_j)-H(T_j))^2] \le1/n$.
This proves (\ref{tech2}) and (\ref{tech3}).
Next, we have almost surely,
\[
\frac{1}{n}\sum_{j=1}^nu_r
\bigl(\H(T_j) \bigr)u_r^{\top} \bigl(
\H(T_j) \bigr)= \frac{1}{n}\sum_{j=1}^nu_r(j/n)u_r^{\top}(j/n).
\]
For a function $h$ defined on $[0,1]$ with Lipschitz constant $L$, we have\vspace*{1pt}
\[
\Biggl|\frac{1}{n}\sum_{j=1}^nh(j/n)-\int
_0^1 h(u)\, \mathrm{d}u \Biggr|\le\frac
{1}{n}
\sum_{j=1}^n\sup_{j-1\le nu \le
j}
\bigl|h(j/n)-h(u)\bigr| \le L/n.
\]\vspace*{1pt}
Since the function $\phi_{k}\phi_{l}$ is Lipschitz with Lipschitz
constant $2\uppi(k+l)$,
we derive the desired bound (\ref{tech1}).
\end{pf}

\begin{details*}[for Example~\ref{ex.model}]
Let $X_1,\dots,X_n$ be independent copies of a random variable $X$
that has distribution function $F_{\theta}$ and density $f_{\theta}$
for some $\theta$ in the open subset $\Theta$ of $\R^q$.
Recall we assumed in Example~\ref{ex.model} that the map $\vt\mapsto
s_{\vt
}=\sqrt{f_{\vt}}$
is continuously differentiable in $L_2$ with derivative $\vt\mapsto
\dot s_{\vt}$
and that the information matrix
$J(\vt)= 4\int\dot s_{\vt}(x) \dot s_{\vt}(x)^{\top}\, \mathrm{d}x$
is invertible for each $\vt$ in~$\Theta$. Thus, we have\vspace*{1pt}
%
\begin{equation}
\label{m-hd} \rho(\tau)=\int \bigl(s_{\theta+\tau}(x)-s_{\theta}(x)-
\tau^{\top} \dot s_{\theta}(x) \bigr)^2\, \mathrm{d}x =
\mathrm{o} \bigl(| \tau|^2 \bigr).
\end{equation}
Recall also that $\ldot_{\theta}= 2\dot s_{\theta}/s_{\theta}$
denotes the score function. By the properties of the densities,
there is a $\delta>0$ and a constant $K$ such that\vspace*{1pt}
%
\begin{equation}
\label{m-dlip} \int|f_{\vt_1}(x)-f_{\vt_2}(x)|\, \mathrm{d}x \leq
K | \vt_1-\vt_2|,\qquad |\vt_1-\theta|<
\delta, | \vt_2-\theta|<\delta.
\end{equation}
As a consequence, we have\vspace*{1pt}
%
\begin{equation}
\label{m-lip} \sup_{x\in\R} \bigl|F_{\vt_1}(x)-F_{\vt_2}(x)\bigr| \leq
K |\vt_1-\vt_2|,\qquad |\vt_1-\theta|<
\delta, | \vt_2-\theta|<\delta.
\end{equation}
%
Let $m=m_n \to\infty$ and $\log(n)m_n^3=\mathrm{o}(n)$.
It suffices to show\vspace*{1pt}
\[
\bigl(-2 \log\el(F_{\hat\theta}) - m_n+q \bigr) / {
\sqrt{2(m_n-q)}} \implies N(0,1).
\]
%
For this, we take $w_n=q_n\circ F_{\theta}$ and
$\hat w_n=q_n\circ F_{\hat\theta}$
with $q_n=(\phi_1,\dots,\phi_{m_n})^{\top}$ and verify (\ref{sa2a})
and (C2) with
$\psi=J(\theta)^{-1/2} \ldot_{\theta}$.
The desired result then follows from (ii) of Corollary~\ref{cor.3}.\vspace*{1pt}

We have $W_n=I_{m_n}=\int\hat w_n\hat w_n^\top \,\mathrm{d}F_{\hat\theta}$
and obtain\vspace*{1pt}
\[
\biggl|\int\hat w_n\hat w_n^{\top} \,
\mathrm{d}F_{\theta} - W_n\biggr| \le2m_n
\int\bigl|f_{\hat\theta}(x)-f_{\theta}(x)\bigr| \,\mathrm{d}x =
\mathrm{o}_p \bigl(m_n^{-1/2} \bigr)
\]
in view of (\ref{m-dlip}) and (\ref{para-eff}).
Thus, (\ref{sa2a}) follows if we verify
%
\begin{equation}
\label{ma2} \biggl|\hat W_n -\bar W_n - \int\hat
w_n\hat w_n^{\top} \,\mathrm{d}F_{\theta}+
W_n \biggr|^2 = \mathrm{o}_p \bigl(m_n^{-1}
\bigr).
\end{equation}
Note that $\psi$ has mean $0$ and identity dispersion matrix under
$F_{\theta}$
and that $A_n \psi$ equals $D_n J(\theta)^{-1} \ldot_{\theta}$,
with $D_n=\int w_n \ldot_{\theta}^{\top}\, \mathrm{d}F_{\theta}$.
Thus, (C2) follows from Remark~\ref{rem.fc},
%
\begin{equation}
\label{m.want} \frac{1}{n}\sum_{j=1}^n
\hat w_n(X_j)-w_n(X_j)+
D_n (\hat\theta -\theta) = \mathrm{o}_p
\bigl(n^{-1/2} \bigr),
\end{equation}
the stochastic expansion (\ref{para-eff}), and the fact that
$|D_n|_o$ is bounded.

We are left to verify (\ref{ma2}) and (\ref{m.want}).
For this, we set
\begin{eqnarray*}
U_{nk}(t)&=& \frac{1}{n}\sum_{j=1}^n
\bigl[\phi_k \bigl(F_{\theta+n^{-1/2}
t}(X_j) \bigr)-
\phi_k \bigl(F_{\theta
}(X_j) \bigr) \bigr],
\\
V_{nkl}(t)&=& \frac{1}{n}\sum_{j=1}^n
\bigl[(\phi_k \phi_l) \bigl(F_{\theta
+n^{-1/2} t}(X_j)
\bigr)-(\phi_k\phi_l) \bigl(F_{\theta}(X_j)
\bigr) \bigr],
\end{eqnarray*}
and note that $D_n^{\top}= (d_1,\dots,d_{m_n})$ with
$d_k = \int\phi_k(F_{\theta}) \ldot_{\theta} \,\mathrm{d}F_{\theta}$.
The statements (\ref{ma2}) and (\ref{m.want}) follow if we show that,
for each finite $C$,
\begin{eqnarray*}
T_{n1}(C)&=& \sup_{|t|\leq C} \sum
_{k=1}^{m_n} \sum_{l=1}^{m_n}
\bigl(V_{nkl}(t)-E \bigl[V_{nkl}(t) \bigr] \bigr)^2
= \mathrm{o}_p \bigl(m_n^{-1} \bigr),
\\
T_{n2}(C)&=&\sup_{|t|\leq C} \sum
_{k=1}^{m_n} \bigl(U_{nk}(t)-E
\bigl[U_{nk}(t) \bigr] \bigr)^2 = \mathrm{o}_p
\bigl(n^{-1} \bigr),
\\
T_{n3}(C)&=& \sup_{|t|\leq C} \sum
_{k=1}^{m_n} \bigl(E \bigl[U_{nk}(t) \bigr]+
n^{-1/2}d_k^{\top}t \bigr)^2 =\mathrm{o}
\bigl(n^{-1} \bigr).
\end{eqnarray*}
The first two statements can be verified using the exponential
inequality given in Lemma~5.2 in Peng and Schick~\cite{PS04}.
This requires the fact that $(\log n)m_n^3 /n \to0$.
%
%

The identity
$f_{\theta+\tau}-f_{\theta} - \ldot_{\theta}^{\top}\tau f_{\theta}
= 2s_{\theta}(s_{\theta+\tau}-s_{\theta} - \dot s_{\theta}^{\top
}\tau)
+ (s_{\theta+\tau}-s_{\theta})^2$
and the definition of $d_k$ yield the formula
\begin{eqnarray*}
\int\phi_k \bigl(F_{\theta}(x) \bigr) \bigl(f_{\theta+\tau}(x)-f_{\theta}(x)
\bigr)\, \mathrm{d}x &=& d_k^{\top}\tau + \int
\phi_k \bigl(F_{\theta}(x) \bigr) \bigl(s_{\theta+\tau}(x)
-s_{\theta}(x) \bigr)^2\, \mathrm{d}x
\\
&&{} + 2\int\phi_k \bigl(F_{\theta}(x) \bigr)
s_{\theta}(x) \bigl(s_{\theta+\tau}(x)-s_{\theta}(x) - \dot
s_{\theta}^{\top}(x) \tau \bigr)\, \mathrm{d}x.
\end{eqnarray*}
In view of this and the fact that $\int\phi_k(F_{\vt})\, \mathrm{d}F_{\vt}=0$
for all $\vt$, we have the identity
\begin{eqnarray*}
E \bigl[U_{nk}(t) \bigr]+d_k^{\top}
tn^{-1/2} &=& \int \bigl(\phi_k \bigl(F_{\theta}(x)
\bigr) -\phi_k \bigl(F_{\theta+n^{-1/2}t}(x) \bigr) \bigr)
\bigl(f_{\theta+n^{-1/2}t}(x)-f_{\theta}(x) \bigr)\, \mathrm{d}x
\\
&& {}- \int\phi_k \bigl(F_{\theta}(x) \bigr)
\bigl(s_{\theta+n^{-1/2}t}(x) -s_{\theta}(x) \bigr)^2\, \mathrm{d}x
\\
&&{} - 2\int\phi_k \bigl(F_{\theta}(x) \bigr)
s_{\theta}(x) \bigl(s_{\theta+n^{-1/2}t}(x)-s_{\theta}(x) -
n^{-1/2} t^{\top} \dot s_{\theta}(x) \bigr) \,\mathrm{d}x.
\end{eqnarray*}
Using (\ref{m-dlip}), (\ref{m-lip}) and the orthonormality of the
the functions $s_{\theta} \phi_k\circ F_{\theta}$, $k=1,2,\dots$\,,
in $L_2$, $T_{n3}(C)$
can be bounded by
\[
\frac{6 \uppi^2 m_n^3 K^4 C^4}{n^2} + 6m_n \biggl(\int \bigl(s_{\theta+n^{-1/2}t}(x)-s_{\theta}(x)
\bigr)^2 \,\mathrm{d}x \biggr)^2 + 12 \sup_{|t|\le C}
\rho \bigl(n^{-1/2}t \bigr). 
\]
The desired statement $T_{n3}(C)=\mathrm{o}(n^{-1})$
now follows from (\ref{m-hd}) and $m_n^3=\mathrm{o}(n)$.
This completes the proof of (\ref{a3}).
\end{details*}

\begin{details*}[for Example~\ref{ex.symm}]
Assume that the distribution function of $X$ is symmetric and continuous.
Then $S=\sign(X)$ and $T=|X|$ are independent,
$S$ has mean zero and variance 1,
and $T$ has a continuous distribution function $H$.
Let $\el$ be defined as in Example~\ref{ex.symm} with as $r=r_n \to
\infty$
and $r_n^3=\mathrm{o}(n)$. It suffices to show that
$(-2 \log\el-(1+r_n))/{\sqrt{2(1+r_n)}}$ is asymptotically standard normal.
This follows from Corollary~\ref{cor.3} if we verify (\ref{a2}) and
(\ref{a3w1}).
These conditions follow from Lemma~\ref{lm.tech} applied
with $S_j=\sign(X_j)$ and $T_j=|X_j|$.
Indeed, in view of the properties of $r_n$,
(\ref{a2}) is a consequence of (\ref{tech1}) and (\ref{a3w1}) of
(\ref{tech3}).
\end{details*}

\begin{details*}[for Example~\ref{ex.indep}]
Assume that $X$ and $Y$ are independent.
Part (a) is an immediate consequence of Corollary~\ref{cor.3a}.
Part (b) follows if we show
$(-2 \log\el(\F,\G)-r_n^2)/{\sqrt{2}r_n}$ is asymptotically
standard normal.
We shall use Corollary~\ref{cor.3} to conclude this.
Here $m_n$ equals $r_n^2$ and thus satisfies $m_n^3=\mathrm{o}(n)$. We shall now verify
(\ref{sa2}) and (\ref{a3w1}).
Let us set
\begin{eqnarray*}
D_{klj}&=&\phi_k \bigl(\F(X_j) \bigr)
\phi_l \bigl(\G(Y_j) \bigr)-\phi_k
\bigl(F(X_j) \bigr)\phi_l \bigl(G(Y_j)
\bigr),
\\
\Phi_{kj}&=& \phi_k \bigl(\F(X_j) \bigr)-
\phi_k \bigl(F(X_j) \bigr) \quad\mbox{and}\quad
\Gamma_{lj}= \phi_l \bigl(\G(Y_j) \bigr)-
\phi_l \bigl(G(Y_j) \bigr).
\end{eqnarray*}
In view of the inequality $|D_{klj}| \le\sqrt{2}|\Phi_{kj}|+\sqrt {2} |\Gamma_{lj}|$,
we obtain with the help of (\ref{tech2}) the bound
\[
\sum_{k=1}^{r_n}\sum
_{l=1}^{r_n}\frac{1}{n}\sum
_{j=1}^nE \bigl[|D_{klj}|^2
\bigr] \le\frac{8\uppi^2 r_n^4}{n}.
\]
From this and $r_n^6=\mathrm{o}(n)$, we conclude (\ref{sa2}).

In view of the identity
$D_{klj}=\phi_k(F(X_j))\Gamma_{lj}+ \phi_l(G(Y_j)) \Phi_{kj}+ \Phi_{kl}\Gamma_{jl}$,
(\ref{a3w1}) follows if we verify
\begin{eqnarray*}
T_{n1}&=&\sum_{k=1}^{r_n}
\sum_{l=1}^{r_n} \Biggl(n^{-1/2}\sum
_{j=1}^n\phi_k
\bigl(F(X_j) \bigr) \Gamma_{lj} \Biggr)^2 =
\mathrm{o}_p(1),
\\
T_{n2}&=& \sum_{k=1}^{r_n}
\sum_{l=1}^{r_n} \Biggl(n^{-1/2}\sum
_{j=1}^n\Phi_{kj}
\phi_l \bigl(G(Y_j) \bigr) \Biggr)^2 =
\mathrm{o}_p(1),
\\
T_{n3}&=&\sum_{k=1}^{r_n}
\sum_{l=1}^{r_n} \Biggl(n^{-1/2}\sum
_{j=1}^n\Phi_{kj}
\Gamma_{lj} \Biggr)^2 =\mathrm{o}_p(1).
\end{eqnarray*}
Applications of (\ref{tech3}) with $S_j=\phi_k(F(X_j))$
yield the bound $E[T_{n1}] \le\uppi^2 r_n^4/n$, and this proves $T_{n1}=\mathrm{o}_p(1)$.
The proof of $T_{n2}=\mathrm{o}_p(1)$ is similar.
To deal with $T_{n3}$, 
we set
\[
H(k,l) = \sum_{j=1}^n
\Phi_{k,j} \Gamma_{lj},\qquad \bar\Phi_k=
\frac{1}{n}\sum_{j=1}^n
\Phi_{kj} 
\quad\mbox{and} \quad\bar\Gamma_{l}=
\frac{1}{n}\sum_{j=1}^n
\Gamma_{lj}. 
\]
Note that $R_j=n\F(X_j)$ is the rank of $X_j$.
Given $Y_1,\dots,Y_n$ and the order statistics $X_{(1)},\dots,X_{(n)}$,
the sum $H(k,l)$ is a simple linear rank statistic with
scores $a(j)= \phi_k(j/n)- \phi_k(F(X_{(j)})$
and coefficients $\G_{lj}$ and consequently has (conditional) mean
$n\bar\Phi_{k} \bar\Gamma_{l}$ and (conditional) variance
\[
\frac{1}{n-1}\sum_{j=1}^n (
\Phi_{kj}-\bar\Phi_k)^2 \sum
_{j=1}^n (\Gamma_{li}-\bar
\Gamma_{l})^2 \le\frac{n}{n-1}\frac{1}{n}\sum
_{j=1}^n\Phi_{kj}^2
\frac{1}{n}\sum_{j=1}^n
\Gamma_{lj}^2.
\]
In view of this bound, we derive the inequality
\[
E[T_{n3}] \le\frac{n}{n-1}\sum_{k=1}^n
E \Biggl[\frac{1}{n}\sum_{j=1}^n
\Phi_{kj}^2 \Biggr] \sum_{l=1}^{r_n}
E \Biggl[\frac{1}{n}\sum_{j=1}^n
\Gamma_{lj}^2 \Biggr]+ n \sum_{k=1}^{r_n}
E \bigl[\bar \Phi_{k}^2 \bigr] \sum
_{l=1}^{r_n} E \bigl[\bar\Gamma_{l}^2
\bigr].
\]
We have
\[
E[\bar\Gamma_k^2)=E \bigl[\bar\Phi_{l}^2
\bigr]= \frac{1}{n} + \Biggl(\frac{1}{n}\sum
_{j=1}^n\phi_k(j/n)
\Biggr)^2 \le \frac{1}{n} + \frac
{2\uppi^2k^2}{n^2}.
\]
Using this and (\ref{tech2}), we obtain $E[T_{n3}]=\mathrm{O}(r_n^6 n^{-2})= \mathrm{o}(1)$
and thus $T_{n3}=\mathrm{o}_p(1)$.
\end{details*}

\section*{Acknowledgements} This work was completed
while Anton Schick was visiting the Department of Statistics
at Texas A\&M University.
He wants to thank the members of the department for their extraordinary
hospitality. Thanks go also to Ingrid Van Keilegom for discussions
and for providing an important reference.
Hanxiang Peng's research was supported in part by
NSF Grant DMS-09-40365.
Anton Schick's research was supported in part by
NSF Grant DMS-09-06551.



%

\printhistory

\end{document}